\documentclass[oneside]{amsart}
\usepackage{amsmath,amssymb,color,amsthm,graphicx,cite}
\usepackage[T1]{fontenc}
\usepackage[utf8]{inputenc}
\usepackage{color,bm}
\usepackage{enumerate} 
\usepackage[backref]{hyperref}
\RequirePackage{luatex85}
\usepackage[all]{xy}
\newtheorem{theorem}{Theorem}[section]
\newtheorem{proposition}[theorem]{Proposition}
\newtheorem{lemma}[theorem]{Lemma}
\newtheorem{corollary}[theorem]{Corollary}
\newtheorem{claim}{Claim}

\theoremstyle{definition}
\newtheorem{definition}{Definition}
\newtheorem{main}{Theorem}

\def\T{\mathbb{T} } 
\def\K{\mathbb{K} } 

\def\F{\mathcal{F} }

\def\R{\mathbb{R} } 
 
\def\Z{\mathbb{Z} } 
\def\nbd{neighborhood } 
 
\def\R{\mathbb{R} } 
 
\def\Sv{\mathop{\mathrm{Sing}}(v)} 
\def\Pv{\mathop{\mathrm{Per}}(v)}

\def\-{\ominus} 
\def\+{\oplus} 
\def\0{\circ}

\author{Tomoo Yokoyama}
\date{\today}
\address{Department of Mathematics, Faculty of Science, Saitama University, Shimo-Okubo 255, Sakura-ku, Saitama-shi, 338-8570 Japan\\}
\email{tyokoyama@rimath.saitama-u.ac.jp}
\thanks{The author was partially supported by JSPS Grant Number 24K06733}

\makeatletter
\@namedef{subjclassname@2020}{%
  \textup{2020} Mathematics Subject Classification}
\makeatother

\usepackage{color,bm}

\title[Coarse chain recurrence \& Morse graphs with finite errors]{Coarse chain recurrence and Morse graphs with finite errors}
\subjclass[2020]{Primary 37B20; Secondary 16W70,37G15,37H20,\\93B03,39A28}

\keywords{Morse decomposition, chain recurrence, bifurcation diagram, filtration, potentials of dynamical systems}

\begin{document}

\begin{abstract}
This paper provides a unified framework connecting dynamical systems with tools from topological data analysis and geometric topology and inspires new interactions among dynamical systems, topology, and nonlinear analysis. To this end, we introduce a one-parameter family of ``chain recurrences'' that generalizes chain recurrence and induces a natural filtration on the underlying metric space of a dynamical system. In particular, the forward directions of the filtrations characterize the level of control required to return to the original position, and the backward directions capture the robustness of the recurrence. The resulting filtrations yield potentials and bifurcation diagrams of dynamical systems that encode the evolution of recurrent sets under bounded total or stepwise perturbations. In addition, we extend Morse graphs to one-parameter families of ``coarse Morse graphs,'' which evolve through vertex collapses reflecting coarse recurrence transitions. These constructions not only refine Conley's decomposition but also reveal singular limit behaviors as the perturbation level vanishes. Furthermore, we establish analogous filtrations for difference equations to bridge the theoretical framework with numerical analysis.
\end{abstract}

\maketitle

\section{Introduction} 

Conley theory states that dynamical systems on compact metric spaces can be decomposed, known as a Morse decomposition, into invariant subsets, each of which is chain recurrent or gradient \cite{conley1978isolated}.
The Morse decomposition induces a directed graph, called the Morse graph, which can capture the gradient behaviors. 
The Morse decompositions for mappings and semiflows are developed for arbitrary metric spaces \cite{franks1988pb,hurley1991chain,hurley1995chain,Rybakowski1987}. 
Moreover, the connected components of the set of chain recurrent points, called chain recurrent components, can be approximated via neighborhoods known as isolating blocks. This approximation algorithm is constructed in \cite{kalies2005algorithmic} and implemented in \cite{arai2009database}.
On the other hand, perturbations sometimes lead to changes in the global structures of dynamical systems.

In the context of control theory for linear systems and nonlinear analysis, as well as differential inclusions in viability theory, the set of states that can be reached using constrained controls is known as the reachable set (cf. \cite{antsaklis2006linear,sontag1998mathematical}) or attainable set (cf. \cite{colonius2000dynamics,jurdjevic1997geometric,aubin1984differential}). 
However, from the perspective of control and prediction, a framework that describes the invariance and sensitivity under perturbations is required. 
In fact, there is a need for a framework that describes properties of escape from attracting basins, the elimination of stagnation, the persistence of recurrent points, and singular limit behaviors where total energy injections go to zero. 
Therefore, to analyze such phenomena and bridge the theory of chain recurrence with reachability, it is necessary to extend the concept of chain recurrence and construct variants of the Morse graph.

Intuitively, this paper establishes a foundation for describing how the global structure of the flow is affected under perturbations of finite energy, by introducing concepts of ``coarse chain recurrence'' (see the next subsection for precise statements).
These concepts of ``coarse chain recurrence'' enable the definition of chain recurrent potentials, which can be utilized to provide a theoretical lower bound for the norm of noise or error required for specific dynamical phenomena to occur in random dynamical systems and stochastic dynamical systems, or in numerical computations of dynamical systems where errors are inevitable.
Furthermore, this natural one-parameter family of ``coarse chain recurrences'' admits a structure where the approximation of chain recurrent components using isolating blocks can be addressed uniformly across the parameter, analogous to how persistent homology addresses a sequence of homology groups.
Moreover, they offer a form of guarantee against the occurrence of such phenomena under insufficient perturbation.
In addition, the concepts of ``coarse chain recurrence'' are intrinsically linked to the solvability of control problems involving quasi-periodicity and recurrence associated with reachable sets. 
In fact, by interpreting the chain errors as admissible inputs, the ``coarse chain recurrence'' of a point $x$ is a necessary condition for $x$ to belong to the reachable set from $x$ itself under the permissible control level.

From the perspective of practical applications and connections to broader mathematical problems, several recent studies provide concrete examples.
For instance, in \cite{pilarczyk2023topological}, the authors utilized rigorous set-oriented numerical methods, Morse decompositions, Morse graphs, and chain recurrent sets to investigate the global dynamics of a discrete neuron model, analyzing the internal structure of attractors via recurrence times. 
Analogously, the concept of recurrent potentials introduced in this paper provides a theoretical framework to conduct a similar analysis of the internal structures of such attractors. 
In addition to its theoretical interest, the approach can also serve as a practical tool for identifying qualitative features in dynamical data. 
Moreover, in \cite{imoto2025filtrations}, an approach incorporating the concept of ``coarse recurrence'' was applied to gradient systems, developing a framework to analyze the sensitivity and robustness of ``attractors'' to perturbations. 
As a practical application to ensemble forecasts of a tropical cyclone, their method reveals spatiotemporal regions of heightened sensitivity, offering new perspectives for topological data analysis.
Furthermore, as suggested in \cite{yokoyama2025coarse}, such analysis of control possibilities is also feasible for problems under the constraint that input-output control can only be performed in the neighborhood of the release point, analogous to ``coarse recurrence''.

\subsection{Statements of main results}

To realize such applications in a rigorous framework, we establish the following main results, which provide filtrations of continuous mappings on metric spaces and one-parameter families of graphs associated with such mappings.

\begin{main}\label{main:01}
For any continuous mapping $f \colon X \to X$ on a metric space $X$, the family $(\mathrm{CR}_{\varepsilon}(f))_{\varepsilon \in \R}$ {\rm(}see Definitions~\ref{def:reachable}, \ref{def:cr}, \ref{def:cr2}, \ref{def:cr02}, and  \ref{def:cr2_neg} below{\rm)} is a filtration with $\mathrm{CR}(f) = \mathrm{CR}_{0}(f)$. 
\end{main}

By using the concept of $\sum$chain recurrence introduced by Easton \cite{easton1978chain}, we obtain the following similar statement, which allows us to describe the injection of finite energy at each step.

\begin{main}\label{main:02}
For any continuous mapping $f \colon X \to X$ on a metric space $X$, the family $(\mathrm{CR}_{\varepsilon}^\Sigma(f))_{\varepsilon \in \R}$ {\rm(}see Definitions~\ref{def:reachable}, \ref{def:cr}, \ref{def:cr2}, \ref{def:cr02}, and  \ref{def:cr2_neg} below{\rm)} is a filtration with $\mathrm{CR}^\Sigma(f) = \mathrm{CR}_{0}^\Sigma(f)$. 
\end{main}

More general statements for these results are described in 
Theorem~\ref{thm:01_gen} for any mappings on metric spaces. 
Although these results hold for mappings and not generally for flows (see Lemma~\ref{lem:no_filtrations_step}), similar results can be derived in the case of semiflows (see Theorem~\ref{th:main03-} and Theorem~\ref{th:main03} below).

We also have the following reductions by using the filtrations $(\mathrm{CR}_{\varepsilon}(f))_{\varepsilon \in \R}$ and $(\mathrm{CR}_{\varepsilon}^\Sigma(f))_{\varepsilon \in \R}$. 

\begin{main}\label{main:02'}
The following statements hold for any continuous mapping $f$ on a metric space and any numbers $\varepsilon_1 < \varepsilon_2 \in \R$: 
\\
{\rm(1)} The Morse graph $G_{(\varepsilon_2,\nu)}(f)$ with $(\varepsilon_2,\nu)$-errors {\rm(}see Definitions~\ref{def:Morse_graph_map02} and \ref{def:MG_error_03} below{\rm)} is partially obtained from the Morse graph $G_{(\varepsilon_1,\nu)}(f)$ by vertex collapses {\rm(i.e.} $G_{(\varepsilon_1,\nu)}(f) \searrow G_{(\varepsilon_2,\nu)}(f)${\rm)} {\rm(}see Definition~\ref{def:ec} below{\rm)}.
\\
{\rm(2)} The Morse graph $G_{(\varepsilon_2,\nu)}^\Sigma(f)$ with $(\varepsilon_2,\nu)$-$\sum$errors 
 is partially obtained from $G_{(\varepsilon_1,\nu)}^\Sigma(f)$ by vertex collapses {\rm(i.e.} $G_{(\varepsilon_1,\nu)}^\Sigma(f) \searrow G_{(\varepsilon_2,\nu)}^\Sigma(f)${\rm)}. 
\end{main}

The previous theorem follows from more general statements described in Theorem~\ref{th:graph_collapse}. 
We observe the following singular limit behaviors, which follow from Theorem~\ref{lem:counter_ex03} and Theorem~\ref{lem:counter_ex02}.

\begin{main}\label{main:counter_ex03}
For any $\varepsilon>0$, there are a homeomorphism $f$ on a metric space and a flow $v$ on a metric space satisfying the following strict inclusions: 
\[
G^\Sigma_{(\varepsilon,0)}(f) \subsetneq \bigcap_{\nu > 0} G^\Sigma_{(\varepsilon,\nu)}(f)
\]
\[
G^\Sigma_{(\varepsilon,0)}(v) \subsetneq \bigcap_{\nu > 0} G^\Sigma_{(\varepsilon,\nu)}(v)
\]
\end{main}

Notice that another type of singular limit of $\mathrm{CR}_{-\varepsilon}^\Sigma$ for mappings exists (see Theorem~\ref{th:counter_ex04}). 

The present paper consists of eleven sections.
In the next section, we recall some concepts of combinatorics and dynamical systems. 
We discuss mappings in \S~3--6 and semiflows in \S~7--9. 
Specifically, the $\varepsilon$-$\ell^p$-chain recurrence, $\varepsilon$-chain recurrence, and $\varepsilon$-$\sum$chain recurrence for mappings are introduced in \S~3, and the chain recurrent potential, ``chain recurrence with negative errors'', and filtrations describing ``coarse recurrence'' are introduced in \S~4.  
In \S~5 (resp. 6), we introduce the Morse (hyper-)graph with finite errors for homeomorphisms (resp. mappings). 
Similarly, in the context of semiflows, we also introduce $\varepsilon$-$\ell^p$-chain recurrence in \S~7, ``chain recurrence with negative errors'' and filtrations describing ``coarse recurrence'' in \S~8, and the Morse (hyper-)graph with finite errors in \S~9. 
In \S~10, we illustrate toy models of escape from attracting basins and elimination of stagnation by controls using finite total energy, persistence of recurrent points, and singular limit behaviors where total energy injections go to zero. 
In the final section, we discuss variants of $\varepsilon$-$\ell^p$-chain recurrence and examine the possible applications.
Additionally, we introduce natural filtrations and ``bifurcation'' diagrams for differential equations and difference equations to bridge the theoretical framework with numerical analysis.

\section{Preliminaries}


We recall basic concepts in set theory, combinatorics,  and dynamical systems.

\subsection{Basic concepts in set theory}

To analyze the various phenomena, we define a filtration as follows.

\begin{definition}\label{def:filtration}
Let $X$ be a set and $\F = \{ F_i \mid i \in I\} \subset 2^X$ a family indexed by a totally ordered set $I$, where $2^X$ is the power set of $X$. The family $\F$ is a {\bf filtration} if $X = \bigcup_{i \in I} F_i$ and $F_{i_1} \subseteq F_{i_2}$ for any pair $i_1 \leq i_2 \in I$.
\end{definition}

\subsection{Basic concepts in combinatorics}

Let $V$ be a set. 
An {\bf (abstract) (simple) directed graph} is a pair of the set $V$ and a subset of $V^2$. 
An {\bf (abstract) directed hyper-graph} is a pair of the set $V$ and a family $H \subseteq V^* \times V^*$, where $V^*$ is the family of non-empty finite subsets of $V$. 
\subsection{Basic concepts in mappings}
Let $f \colon X \to X$ be a mapping on a metric space $X$, which need not be continuous. 
For a point $x \in X$, we denote by $O^+(x)$ the positive orbit (i.e. $O^+(x) := \{ f^n(x) \mid n \in \Z_{> 0} \}$). 
A point $x \in X$ is {\bf fixed} if $f(x) = x$, and is {\bf periodic} if there is a positive integer $n$ with $f^n(x)= x$.
The $\omega$-limit set of a point $x$ is defined as follows: 
\[
\begin{split}
\omega(x) := \bigcap_{n\in \mathbb{Z}}\overline{\{f^m(x) \mid m > n\}} = \bigcap_{n\in \mathbb{Z}} \overline{f(\Z_{>n},x)} &
\end{split}
\]

\begin{definition}
For any positive integer $n \in \Z_{>0}$ and any non-negative number $\varepsilon \geq 0$, a sequence $( x_i )_{i=0}^{n}$ of points $x_i \in X$ is an {\bf $\bm{\varepsilon}$-chain} if there is a sequence $( \varepsilon_i )_{i=0}^{n-1}$ of non-negative numbers with $\varepsilon \geq \max_{i \in \{0, \ldots , n-1\}} \{ \varepsilon_i \}$ such that $d(f(x_i), x_{i+1}) \leq \varepsilon_i$ for any $i \in \{ 0,1, \ldots , n-1 \}$. 
\end{definition}

For any $\varepsilon \geq 0$, define a binary relation $\sim_\varepsilon$ as follows: 
\[
x \sim_\varepsilon y \text{ if there is an } \varepsilon \text{-chain from }x \text{ to } y  
\]

Define a binary relation $\sim_{\mathrm{CR}}$ by $x \sim_{\mathrm{CR}} y$ if $x \sim_{\varepsilon} y$ for any $\varepsilon > 0$. 
Define a binary relation $\approx_{\mathrm{CR}}$ by $x \approx_{\mathrm{CR}} y$ if $x \sim_{\mathrm{CR}} y$ and $y \sim_{\mathrm{CR}} x$.
Notice that the binary relation $\approx_{\mathrm{CR}}$ is an equivalence relation. 
We have the following variant of recurrence. 

\begin{definition}
A point $x \in X$ is {\bf chain recurrent} if $x \sim_\varepsilon x$ for any $\varepsilon > 0$. 
The set of chain recurrent points is called the {\bf chain recurrent set} of $f$ and is denoted by $\bm{\mathrm{CR}(f)}$. 
\end{definition}

We call an equivalence class of the chain recurrent set $\mathrm{CR}(f)$ with respect to $\approx_{\mathrm{CR}}$ a {\bf ($\bm{0}$-)chain recurrent component} of $\mathrm{CR}(f)$.

\subsubsection{Variant of chain recurrence}

We recall the following strong $\varepsilon$-chain introduced by Robert Easton \cite{easton1978chain}, 
which is also called a summable chain or $\sum$chain in \cite{de2025streams} by Roberto De Leo and James Yorke. 

\begin{definition}
For any positive integer $n \in \Z_{>0}$ and any non-negative number $\varepsilon \geq 0$, a sequence $( x_i )_{i=0}^n$ of points $x_i \in X$ is an {\bf $\bm{\varepsilon}$-$\sum$chain} {\rm(}or {\bf strong $\bm{\varepsilon}$-chain}{\rm)} if there is a sequence $( \varepsilon_i )_{i=0}^{n-1}$ of non-negative numbers with $\varepsilon \geq \sum_{i=0}^{n-1} \varepsilon_i$ such that $d(f(x_i), x_{i+1}) \leq \varepsilon_i$ for any $i \in \{ 0,1, \ldots , n-1 \}$. 
\end{definition}

For any $\varepsilon \geq 0$, define a binary relation $\sim_{\varepsilon_\Sigma}$ as follows: 
\[
x \sim_{\varepsilon_\Sigma} y \text{ if there is an  } \varepsilon\text{-}{\textstyle \sum} \text{chain from }x \text{ to } y  
\]
Define a binary relation $\sim_{\mathrm{CR}^\Sigma}$ by $x \sim_{\mathrm{CR}^\Sigma} y$ if $x \sim_{\varepsilon_\Sigma} y$ for any $\varepsilon > 0$. 
Define a binary relation $\approx_{\mathrm{CR}^\Sigma}$ by $x \approx_{\mathrm{CR}^\Sigma} y$ if $x \sim_{\mathrm{CR}^\Sigma} y$ and $y \sim_{\mathrm{CR}^\Sigma} x$.
Notice that the binary relation $\approx_{\mathrm{CR}^\Sigma}$ is an equivalence relation. 
We have the following variant of recurrence. 

\begin{definition}
A point $x \in X$ is {\bf $\sum$chain recurrent} if $x \sim_{\mathrm{CR}^\Sigma} x$. 
Denote by $\bm{\mathrm{CR}^\Sigma(f)}$ the set of $\sum$chain recurrent points. 
\end{definition}

We call an equivalence class of the chain recurrent set $\mathrm{CR}^\Sigma(f)$ with respect to $\approx_{\mathrm{CR}^\Sigma}$ a {\bf $\bm{0}$(-$\sum)$chain recurrent component} of $\mathrm{CR}^\Sigma(f)$.

\subsection{Basic concepts for homeomorphisms}

Let $f \colon X \to X$ be a homeomorphism on a topological space $X$.
For a point $x \in X$, we denote by $O(x)$ the orbit of $x$ by $f$ and $O^-(x)$ the negative orbit (i.e. $O^-(x) := \{ f^n(x) \mid n \in \Z_{< 0} \}$).
The $\alpha$-limit set of a point $x$ is defined as follows: 
\[
\alpha(x) := \bigcap_{n \in \mathbb{Z}}\overline{\{f^t(x) \mid t < n \}} = \bigcap_{n \in \mathbb{Z}} \overline{f(\Z_{< n},x)}
\]

A point $x$ of $X$ is {\bf recurrent} (cf. \cite{marzougui2002area}) if $x \in \omega(x) \cup \alpha(x)$.
Notice that some authors refer to the property $x \in \omega(x) \cap \alpha(x)$ as recurrence (cf. \cite{nikolaev2013foliations}).

\subsubsection{Morse graph for homeomorphisms}\label{sec:Morse_graph}

Consider a set $\mathcal{X} = \{ X_i \}_{i \in \Lambda}$ of disjoint subsets. 
For any $j \neq k \in \Lambda$, define subsets $D_{j,k}'$ and $D_{j,k}$ as follows: 
\[
\begin{split}
D'_{j,k} &:= \left\{ x \in X - \bigsqcup_{i \in \Lambda} X_i \middle| \, \alpha(x) \cap  X_j \neq \emptyset,  \omega(x) \cap  X_k \neq \emptyset \right\}
\\
D_{j,k} &:= \left\{ x \in X \middle| \, \alpha(x) \cap  X_j \neq \emptyset,  \omega(x) \cap  X_k \neq \emptyset \right\}
\end{split}
\]
where the symbol $\bigsqcup$ denotes a disjoint union and the symbol $A - B$ is used instead of the set difference $A \setminus B$ when $B \subseteq A$.
A directed graph $(V, D')$ is called the {\bf Morse graph} of $\mathcal{X}$ if $V = \{ X_i \mid i  \in \Lambda \}$ and $D' = \{ (X_j, X_k) \mid D'_{j,k} \neq \emptyset, j \neq k \in \Lambda \}$. 
Then the Morse graph $(V, D')$ is denoted by $\bm{G'_{\mathcal{X}}(f)}$. 
Similarly, a directed graph $(V, D)$ is called the {\bf Morse graph} of $\mathcal{X}$ equipped with the  {\bf persistent directed edge set} $D$ if $V = \{ X_i \mid i  \in \Lambda \}$ and $D = \{ (X_j, X_k) \mid D_{j,k} \neq \emptyset, j \neq k \in \Lambda \}$. 
Then the Morse graph $(V, D)$ is denoted by $\bm{G_{\mathcal{X}}(f)}$. 

\begin{definition}\label{def:Morse_graph_map_00}
The graph $G'_{\mathcal{X}}(f)$ (resp. $G_{\mathcal{X}}(f)$) is the {\bf Morse graph} (resp. {\bf Morse graph} with persistent directed edges) of $f$ if $\mathcal{X}$ is the set of $0$-chain recurrent components of $\mathrm{CR}(f)$. 
Then we denote the Morse graph of $f$ as $\bm{G'(f)}$ (resp. $\bm{G(f)}$). 
\end{definition}

We have the following observation. 

\begin{lemma}\label{lem:edge_corr}
Let $f \colon X \to X$ be a homeomorphism on a compact metric space $X$ and denote by $\mathcal{X}= \{ X_i \}_{i \in \Lambda}$ the set of chain recurrent components of $\mathrm{CR}(f)$. 
Then $G'(f) = G(f)$.
\end{lemma}

\begin{proof}
Since the chain recurrent set $\mathrm{CR}(f)$ is compact $f$-invariant, for any equivalence class $X_i$ of $\mathrm{CR}(f)$ and any $y \in X_i$, we have $\alpha(y) \cup \omega(y) \subset X_i$. 
For any $j \neq k \in \Lambda$, we have $D'_{j,k} \cap \mathrm{CR}(f) = \emptyset$ and so $D_{j,k} = D'_{j,k}$. 
\end{proof}

Although $G'(f)$ is the original Morse graph and $G(f)$ is a variant of $G'(f)$, we introduce $G(f)$ because the persistence of directed edges is required for constructing a filtration (see \S~\ref{sec:non-increasing} for non-persistence of directed edges of $G'(f)$).

\subsection{Basic concepts for semiflows}

A {\bf flow} is a continuous $\R$-action on a topological space.
A {\bf semiflow} is a continuous $\R_{\geq 0}$-action on a topological space.

Let $v$ be a continuous action $v \colon \K \times X \to X$, where $\K$ is either $\R$ or $\R_{\geq 0}$. 
Notice that if $\K = \R$ (resp. $\K = \R_{\geq 0}$) then $v$ is a flow (resp. semiflow). 
For $t \in \K$, define $v^t : X \to X$ by $v^t := v(t, \cdot )$.

A point $x \in X$ is {\bf singular} if $x = v^t(x)$ for any $t \in \K$. 
For a point $x \in X$, we denote by $O(x)$ the orbit of $x$ by $v$ (i.e. $O(x) := v(\K,x)$), and the positive orbit (i.e. $O^+(x) := \{ v^t(x) \mid t > 0 \} = v(\R_{>0},x)$). 
The $\omega$-limit set of a point $x$ is defined as follows: 

\[
\begin{split}
\omega(x) := \bigcap_{r \in \K}\overline{\{v^t(x) \mid t > r\}} = \bigcap_{r \in \K} \overline{v(\R_{>r},x)} &
\end{split}
\]

\subsubsection{Chain recurrence and related concepts for semiflows}\label{cr_semiflow001}

We recall some concepts to define the chain recurrence. 

\begin{definition}
For any positive integer $n \in \Z_{>0}$ and any non-negative numbers $\varepsilon \geq 0$ and $T>0$, a sequence $( x_i )_{i=0}^n$ of points $x_i \in X$ is an {\bf $\bm{(\varepsilon,T)}$-chain} if there are sequences $(\varepsilon_i)_{i=0}^{n-1}$, $(t_i)_{i=0}^{n-1}$ of non-negative numbers with $\varepsilon \geq \max_{i \in \{ 0, \ldots , n-1\}} \{ \varepsilon_i \}$ and $\min_i \{ t_i \} \geq T$  such that $d(v^{t_i}(x_i), x_{i+1}) \leq \varepsilon_i$ for any $i \in \{ 0,1, \ldots , n-1 \}$. 
\end{definition}

For any $\varepsilon \geq 0$ and $T>0$, define a binary relation $\sim_{(\varepsilon,T)}$ as follows: 
\[
x \sim_{(\varepsilon,T)} y \text{ if there is an } (\varepsilon,T) \text{-chain from }x \text{ to } y  
\]
We have the following chain recurrence. 

\begin{definition}
A point $x \in X$ is {\bf chain recurrent} if $x \sim_{(\varepsilon,T)} x$ for any $\varepsilon > 0$ and $T>0$. 
Denote by $\bm{\mathrm{CR}(v)}$ the set of chain recurrent points of the semiflow $v$. 
\end{definition}

Notice that Morse graphs $G_{\mathcal{X}}(v)$ and $G(v)$ for flows are defined in the same way as for homeomorphism in \S~\ref{sec:Morse_graph}. 


\subsubsection{Variant of the chain recurrence}
We define the $\sum$chain recurrence as follows. 

\begin{definition}
For any positive integer $n \in \Z_{>0}$ and any non-negative number $\varepsilon \geq 0$, a sequence $( x_i )_{i=0}^n$ of points $x_i \in X$ is a {\bf $(\varepsilon,T)$-$\sum$chain} if there are sequences $(\varepsilon_i)_{i=0}^{n-1}$, $(t_i)_{i=0}^{n-1}$ of non-negative numbers with $\varepsilon \geq \sum_{i = 0}^{n-1} \varepsilon_i$ and $\min_i \{ t_i \} \geq T$  such that $d(v^{t_i}(x_i), x_{i+1}) \leq \varepsilon_i$ for any $i \in \{ 0,1, \ldots , n-1 \}$. 
\end{definition}

For any $\varepsilon \geq 0$ and $T>0$, define a binary relation $\sim_{\varepsilon_\Sigma}$ as follows: 
\[
x \sim_{(\varepsilon_\Sigma,T)} y \text{ if there is an } (\varepsilon,T)\text{-}{\textstyle \sum}\text{chain from }x \text{ to } y  
\]

\begin{definition}
A point $x \in X$ is {\bf $\sum$chain recurrent} if $x \sim_{(\varepsilon_\Sigma,T)} x$ for any $\varepsilon > 0$ and $T>0$. 
Denote by $\bm{\mathrm{CR}^\Sigma(v)}$ the set of $\sum$chain recurrent points of the semiflow $v$. 
\end{definition}

Note that a variant of chain recurrence for (semi-)flows is discussed in \S~\ref{rem:f_rem}. 

\subsubsection{Basic concepts for flows}\label{sec:flow}

Suppose that $v$ is a flow in this \S~\ref{sec:flow}. 
For a point $x \in X$, we denote by $O^-(x)$ the negative orbit (i.e. $O^-(x) := \{ v^t(x) \mid t < 0 \} = v(\R_{<0},x)$).
The $\alpha$-limit set of a point $x$ is defined as follows: 
\[
\alpha(x) := \bigcap_{r \in \mathbb{R}}\overline{\{v^t(x) \mid t < r\}} = \bigcap_{r \in \mathbb{R}} \overline{v(\R_{<r},x)}
\] 
A point $x$ of $X$ is {\bf recurrent} if $x \in \omega(x) \cup \alpha(x)$.

\section{Coarse chain recurrence for mappings} 

In this section, we introduce $\varepsilon$-chain recurrence and $\varepsilon$-$\sum$chain recurrence, as well as a more general concept. 
Let $(X,d)$ be a metric space and $f \colon X \to X$ a mapping.

\subsection{$\varepsilon$-$\ell^p$-chains and $\ell^p$-chain recurrence}

Since the error for an $\varepsilon$-chain can be seen as being evaluated by the $\ell^\infty$ norm and the error for an $\varepsilon$-$\sum$chain is evaluated by the $\ell^1$ norm\footnote{This observation is suggested by Yusuke Imoto}, we define the $\varepsilon$-$\ell^p$-chains to unify these concepts as follows.

\begin{definition}
For any $p \in [1,\infty]$, any positive integer $n \in \Z_{>0}$, and any non-negative number $\varepsilon \geq 0$, a sequence $( x_i )_{i=0}^n$ of points $x_i \in X$ is an {\bf $\bm{\varepsilon}$-$\bm{\ell^p}$-chain} from a point $x \in X$ to a point $y \in X$ if $x_0 = x$ and $x_n = y$ and  there is a sequence $( \varepsilon_i )_{i=0}^{n-1}$ of non-negative numbers with $\varepsilon \geq \| ( \varepsilon_i )_{i=0}^{n-1} \|_{p}$ such that $d(f(x_i), x_{i+1}) = \varepsilon_i$ for any $i \in \{ 0,1, \ldots , n-1 \}$, where 
\[
\| ( \varepsilon_i )_{i=0}^{n-1} \|_{p} := 
\begin{cases}
\left( \vert \varepsilon_0 \vert^p  + \vert \varepsilon_1 \vert^p + \cdots + \vert \varepsilon_{n-1}\vert^p \right)^{1/p} &  (p \in [1,\infty) \\
\max\{ \varepsilon_0, \varepsilon_1, \ldots , \varepsilon_{n-1} \} &  (p = \infty) 
\end{cases}
\]
is the $\ell^p$-norm. 
\end{definition}

We have the following binary relation. 

\begin{definition}
For any $p \in [1,\infty]$ and any $\varepsilon \geq 0$, we define a binary relation $\sim_{\varepsilon +}$ as follows: 
\[
x \sim^{\ell^p}_{\varepsilon} y \text{ if } \text{there is an } \varepsilon \text{-}\ell^p\text{-}\text{chain from }x \text{ to } y\]
\[
x \sim^{\ell^p}_{\varepsilon +} y \text{ if } x \sim^{\ell^p}_{\varepsilon'} y \text{ for any } \varepsilon' > \varepsilon 
\]
\end{definition}

We define $\varepsilon$-reachability (or $\varepsilon$-attainability) as follows. 

\begin{definition}\label{def:reachable}
For any $p \in [1,\infty]$ and any $\varepsilon \geq 0$, define an {\bf $\bm{\varepsilon}$-$\ell^p$-reachable} (or {\bf $\bm{\varepsilon}$-$\ell^p$-attainable}) {\bf set} $[x]^{\ell^p}_\varepsilon$ of $x$, a weakly {\bf $\bm{\varepsilon}$-$\ell^p$-reachable} (or {\bf $\bm{\varepsilon}$-$\ell^p$-attainable}) {\bf set} $[x]^{\ell^p}_{\varepsilon+}$ of $x$, and  a weakly {\bf $\bm{-\varepsilon}$-$\ell^p$-reachable} (or {\bf $\bm{-\varepsilon}$-$\ell^p$-attainable}) {\bf set} $[x]^{\ell^p}_{-\varepsilon+}$ of $x$  as follows: 
\[
[x]^{\ell^p}_\varepsilon := \{ y \in X \mid x \sim^{\ell^p}_\varepsilon y \} = \{ y \in X \mid \text{There is an } \varepsilon\text{-}\ell^p\text{-chain from }x \text{ to } y \}	
\]
\[
[x]^{\ell^p}_{\varepsilon+} := \bigcap_{\varepsilon' > \varepsilon} [x]^{\ell^p}_{\varepsilon'} = \bigcap_{\varepsilon' > \varepsilon}  \{ y \in X \mid x \sim^{\ell^p}_ {\varepsilon'} y \} = \{ y \in X \mid x \sim^{\ell^p}_{\varepsilon+} y \}
\]
\[
[x]^{\ell^p}_{-\varepsilon+} := \bigcap_{\varepsilon' \in [0,\varepsilon]} \bigcap_{z \in [x]^{\ell^p}_{\varepsilon'}} [z]^{\ell^p}_{\varepsilon'+} = \bigcap_{\varepsilon' \in [0,\varepsilon]} \bigcap_{z \in [x]^{\ell^p}_{\varepsilon'}} \bigcap_{\varepsilon'' > \varepsilon'} [z]^{\ell^p}_{\varepsilon''} 
\]
\end{definition}

Intuitively, the subset $[x]^{\ell^p}_{-\varepsilon+}$ is the set of points that can return using controls of magnitude almost $\varepsilon$, under any perturbation of magnitude $\varepsilon$. 
Note that we define $[x]^{\ell^p}_{-\varepsilon+}$ here to allow for comparison with other notation, although it will be discussed in the next section.
We have the characterizations of chain recurrence and $\sum$chain recurrence as follows. 

\begin{lemma}\label{lem:ell_infty}
The following conditions are equivalent for any point $x \in X$: 
\\
{\rm(1)} The point $x$ is chain recurrent. 
\\
{\rm(2)} For any $\varepsilon' > 0$, there is an $\varepsilon'$-$\sum$chain from $x$ to $x$. 	
\\
{\rm(3)} $x \sim^{\ell^\infty}_{0 +} x$. 
\\
{\rm(4)} $x \in [x]^{\ell^\infty}_{0 +}$. 
\end{lemma}

\begin{lemma}\label{lem:ell_1}
The following conditions are equivalent for any point $x \in X$: 
\\
{\rm(1)} The point $x$ is $\sum$chain recurrent. 
\\
{\rm(2)} For any $\varepsilon' > 0$, there is an $\varepsilon'$-chain from $x$ to $x$. 	
\\
{\rm(3)} $x \sim^{\ell^1}_{0 +} x$. 
\\
{\rm(4)} $x \in [x]^{\ell^1}_{0 +}$. 
\end{lemma}

By the previous lemma, we introduce the following variant of chain recurrence. 

\begin{definition}
For any $p \in [1,\infty]$, a point $x \in X$ is {\bf $\bm{\ell^p}$-chain recurrent} if $x \in [x]^{\ell^p}_{0 +}$. 
Denote by $\bm{\mathrm{CR}^{\ell^p}(f)}$ the set of $\ell^p$-chain recurrent points. 
\end{definition}

\subsubsection{Variants of chain recurrence}

We have the following observation.

\begin{lemma}
For any $p \in [1,\infty]$, a point $x \in X$ is periodic if and only if $x \sim^{\ell^p}_{0} x$. 
\end{lemma}

By the previous lemma, we have the following variant of periodicity. 

\begin{definition}
For any $p \in [1,\infty]$, a point $x \in X$ is {\bf $\bm{\varepsilon}$-$\ell^p$-periodic} if $x \sim^{\ell^p}_{\varepsilon} x$ {\rm(i.e.} $x \in [x]^{\ell^p}_{\varepsilon}${\rm)}. 
\end{definition}

Denote by $\bm{\mathop{\mathrm{Per}^{\ell^p}_{\varepsilon}}(f)}$ the set of $\varepsilon$-$\ell^p$-periodic points. 
We have the following observation. 

\begin{lemma}\label{lem:equality0}
For any $p \in [1,\infty]$, the family $(\mathop{\mathrm{Per}^{\ell^p}_{\varepsilon}}(f))_{\varepsilon \geq 0}$ is a filtration. 
In fact, the following properties hold: 
\\
{\rm(1)} For any non-negative numbers $\varepsilon_1 < \varepsilon_2$, we have $\mathop{\mathrm{Per}^{\ell^p}_{\varepsilon_1}}(f) \subseteq \mathop{\mathrm{Per}^{\ell^p}_{\varepsilon_2}}(f)$. 
\\
{\rm(2)} $X = \bigcup_{\varepsilon \geq 0} \mathop{\mathrm{Per}^{\ell^p}_{\varepsilon}}(f)$. 
\end{lemma}

\subsection{Coarse chain recurrence}

We have the following observations. 

\begin{lemma}\label{lem:001}
For any $p \in [1,\infty]$ and any non-negative numbers $\varepsilon_1 \leq \varepsilon_2$, we have $[x]^{\ell^p}_{\varepsilon_1} \subseteq [x]^{\ell^p}_{\varepsilon_2}$ and $[x]^{\ell^p}_{\varepsilon_1 +} \subseteq [x]^{\ell^p}_{\varepsilon_2 +}$. 
\end{lemma}

\begin{lemma}\label{lem:002}
We have the following statements for any $p \in [1,\infty]$ and any non-negative numbers $\varepsilon_1 < \varepsilon_2 < \varepsilon_3$ and for any $x \in X$: 
\[
[x]^{\ell^p}_{\varepsilon_1 +} \subseteq [x]^{\ell^p}_{\varepsilon_2} \subseteq [x]^{\ell^p}_{\varepsilon_2 +} \subseteq [x]^{\ell^p}_{\varepsilon_3}
\]
\end{lemma}

\begin{proof}
Lemma~\ref{lem:001} implies $[x]^{\ell^p}_{\varepsilon_2} \subseteq [x]^{\ell^p}_{\varepsilon_2 +}$.
By definitions, we have $[x]^{\ell^p}_{\varepsilon_2 +} = \{ y \in X \mid x \sim^{\ell^p}_{\varepsilon_2 +} y \} \subseteq \{ y \in X \mid x \sim^{\ell^p}_{\varepsilon_3} y \} = [x]^{\ell^p}_{\varepsilon_3 +}$. 
The same argument implies $[x]^{\ell^p}_{\varepsilon_1 +} \subseteq [x]^{\ell^p}_{\varepsilon_2}$. 
\end{proof}

By Lemma~\ref{lem:ell_infty} and Lemma~\ref{lem:ell_1}, we introduce the following variant of chain recurrence to capture behaviors under finite perturbations. 

\begin{definition}\label{def:cr}
For any $p \in [1,\infty]$ and any number $\varepsilon \geq 0$, a point $x \in X$ is {\bf $\bm{\varepsilon}$-$\bm{\ell^p}$-chain recurrent} if $x \in [x]^{\ell^p}_{\varepsilon +}$. 
Denote by $\bm{\mathrm{CR}^{\ell^p}_{\varepsilon}(f)}$ the set of $\varepsilon$-$\ell^p$-chain recurrent points. 
\end{definition}

\begin{definition}\label{def:cr2}
The subset $\mathrm{CR}^{\ell^\infty}_{\varepsilon}(f)$ is denoted by $\bm{\mathrm{CR}_{\varepsilon}(f)}$ and its points are called {\bf $\bm{\varepsilon}$-chain recurrent points}.
Similarly, the subset $\mathrm{CR}^{\ell^1}_{\varepsilon}(f)$ is denoted by $\bm{\mathrm{CR}^\Sigma_{\varepsilon}(f)}$ and its points are called {\bf $\bm{\varepsilon}$-$\bm{\sum}$chain recurrent points}.
\end{definition}

Notice that a point is $0$-$\ell^\infty$-chain recurrent if and only if it is chain recurrent. 
In addition, a point is $0$-$\ell^1$-chain recurrent if and only if it is strong chain recurrent in the sense of Easton \cite{easton1978chain} (i.e. $\sum$chain recurrent). 
%
Moreover, note that $\mathrm{CR}^\Sigma(f) = \mathrm{CR}_{0}^{\ell^1}(f) \subseteq \mathrm{CR}_{0}^{\ell^\infty}(f) = \mathrm{CR}(f)$. 
By definition of $\varepsilon$-$\ell^p$-chain recurrence, we have the following equivalence. 

\begin{lemma}
The following conditions are equivalent for any $p \in [1,\infty]$, any number $\varepsilon \geq 0$,  and any point $x \in X$: 
\\
{\rm(1)} The point $x$ is $\varepsilon$-$\ell^p$-chain recurrent. 
\\
{\rm(2)} $x \sim^{\ell^p}_{\varepsilon +} x$. 	
\\
{\rm(3)} For any $\varepsilon' > \varepsilon$, there is an $\varepsilon'$-$\ell^p$-chain from $x$ to $x$. 
\end{lemma}

Notice that $\mathrm{CR}^{\ell^p}(f) = \mathrm{CR}^{\ell^p}_{0}(f) \subseteq \mathrm{CR}^{\ell^p}_{\varepsilon}(f)$ for any $p \in [1,\infty]$ and any $\varepsilon \geq 0$. 
We define the $\varepsilon$-chain recurrent components to define the Morse (hyper-)graphs with $\varepsilon$-$\ell^p$-errors. 

\begin{definition}\label{def:crc_map}
Define a binary relation $\approx_{\mathrm{cr}^{\ell^p}_{\varepsilon+}}$ on $\mathrm{CR}^{\ell^p}_{\varepsilon}(f)$ by $x \approx_{\mathrm{cr}^{\ell^p}_{\varepsilon+}} y$ if $x \sim^{\ell^p}_{\varepsilon+} y$ and $y \sim^{\ell^p}_{\varepsilon+} x$.
The transitive closure of $\approx_{\mathrm{cr}^{\ell^p}_{\varepsilon +}}$ on $\mathrm{CR}^{\ell^p}_{\varepsilon}(f)$ is an equivalence relation on $\mathrm{CR}^{\ell^p}_{\varepsilon}(f)$, denote by $\approx_{\mathrm{CR}^{\ell^p}_{\varepsilon +}}$. 
We call an equivalence class of the chain recurrent set $\mathrm{CR}^{\ell^p}_{\varepsilon}(f)$ with respect to $\approx_{\mathrm{CR}^{\ell^p}_{\varepsilon+}}$ an {\bf $\bm{\varepsilon}$-{\rm(}$\bm{\ell^p}$-{\rm)}chain recurrent component} of $\mathrm{CR}^{\ell^p}_{\varepsilon}(f)$. 
\end{definition}

%
In addition, we have the following observation. 

\begin{lemma}\label{lem:altanative}
We have the following statement for any $p \in [1,\infty]$ and any non-negative numbers $\varepsilon_1 < \varepsilon_2 < \varepsilon_3$:
\[
\mathrm{CR}^{\ell^p}_{\varepsilon_1 }(f) \subseteq \mathop{\mathrm{Per}^{\ell^p}_{\varepsilon_2}}(f) \subseteq \mathrm{CR}^{\ell^p}_{\varepsilon_2}(f) \subseteq \mathop{\mathrm{Per}^{\ell^p}_{\varepsilon_3}}(f)
\]
\end{lemma}

\begin{proof}
Fix a point $x \in \mathrm{CR}^{\ell^p}_{\varepsilon_1 }(f)$. 
Lemma~\ref{lem:002} implies $x \in [x]^{\ell^p}_{\varepsilon_1 +} \subseteq [x]^{\ell^p}_{\varepsilon_2}$. 
This means that $x \in \mathop{\mathrm{Per}^{\ell^p}_{\varepsilon_2}}(f)$. 

Fix a point $x \in \mathop{\mathrm{Per}^{\ell^p}_{\varepsilon_2}}(f)$. 
Lemma~\ref{lem:002} implies $x \in [x]^{\ell^p}_{\varepsilon_2} \subseteq [x]^{\ell^p}_{\varepsilon_2 +}$ and so $x \in \mathrm{CR}^{\ell^p}_{\varepsilon_2}(f)$.
\end{proof}

By Lemma~\ref{lem:altanative}, notice that filtrations $(\mathop{\mathrm{Per}^{\ell^p}_{\varepsilon}}(f))_{\varepsilon \geq 0}$ and $(\mathrm{CR}^{\ell^p}_{\varepsilon}(f))_{\varepsilon \geq 0}$ are alternately inclusive.
Although the differences between the filtrations may exhibit theoretically intriguing properties, they are virtually indistinguishable from the perspective of data analysis that deals only with finite structures.
Therefore, from now on, we only focus on the filtration $(\mathrm{CR}^{\ell^p}_{\varepsilon}(f))_{\varepsilon \geq 0}$ in this paper. 

We obtain the following observation. 

\begin{lemma}\label{lem:equality}
For any $p \in [1,\infty]$, the family $(\mathrm{CR}^{\ell^p}_{\varepsilon}(f))_{\varepsilon \geq 0}$ is a filtration. 
In fact, the following properties hold: 
\\
{\rm(1)} For any non-negative numbers $\varepsilon_1 < \varepsilon_2$, we have $\mathrm{CR}^{\ell^p}_{\varepsilon_1 }(f) \subseteq \mathrm{CR}^{\ell^p}_{\varepsilon_2}(f)$. 
\\
{\rm(2)} $X = \bigcup_{\varepsilon \geq 0} \mathrm{CR}^{\ell^p}_{\varepsilon}(f)$. 
\end{lemma}

When a value \(\varepsilon_0\) satisfies $X = \mathrm{CR}^{\ell^p}_{\varepsilon_0}(f)$, it seems that the dynamical system mixes under control or noise of the norm \(\varepsilon_0\). 
Therefore, this filtration $(\mathrm{CR}^{\ell^p}_{\varepsilon}(f))_{\varepsilon \geq 0}$ can also be considered as an indicator of the circulation of the dynamical system.
Hence we define the following concept. 

\begin{definition}
We call the inferior $\inf \{ \varepsilon \geq 0 \mid \mathrm{CR}^{\ell^p}_{\varepsilon}(f) = X \}$ the {\bf $\ell^p$-circulation cost} of $f$. 
\end{definition}

We have the following observation.

\begin{proposition}
For any $p \in [1,\infty]$, we have the following equality 
$\mathrm{CR}^{\ell^p}(f) = \mathrm{CR}_{0}^{\ell^p}(f) = \bigcap_{\varepsilon > 0} \mathrm{CR}^{\ell^p}_{\varepsilon}(f)$. 
\end{proposition}

\begin{proof}
By Lemma~\ref{lem:equality}, we have $\mathrm{CR}_{0}^{\ell^p}(f) \subseteq \mathrm{CR}_{\varepsilon}^{\ell^p}(f)$ for any $\varepsilon > 0$. 
Fix any point $x \in \bigcap_{\varepsilon > 0} \mathrm{CR}_{\varepsilon}^{\ell^p}(f)$. 
For any $\varepsilon > 0$, we have $x \sim^{\ell^p}_{\varepsilon +} x$ and so $x \sim^{\ell^p}_{\varepsilon '} x$ for any $\varepsilon' > \varepsilon$. 
This means that $x \sim^{\ell^p}_{\varepsilon} x$ for any $\varepsilon > 0$, and so that $x \sim^{\ell^p}_{0 +} x$. 
Therefore, we obtain $x \in \mathrm{CR}_{0}^{\ell^p}(f)$. 
\end{proof}

\section{Non-gradient property for mappings}

To extend the filtration $(\mathrm{CR}^{\ell^p}_{\varepsilon}(f))_{\varepsilon \geq 0}$, we introduce a non-gradient property as follows. 
Let $(X,d)$ be a metric space and $f \colon X \to X$ a mapping. 

\subsection{$\varepsilon$-non-gradient property}

We define the following binary relation to formulate ``chain recurrence with negative errors''. 

\begin{definition}\label{def:neg_grad_map_001}
For any $p \in [1,\infty]$ and any $\varepsilon \geq 0$, we define a binary relation $\sim^{\ell^p}_{-\varepsilon +}$ as follows: 
\[
x \sim^{\ell^p}_{-\varepsilon +} y \text{ if } z \sim^{\ell^p}_{\varepsilon' +} y \text{ for any } \varepsilon' \in [0,\varepsilon] \text{ and for any } z \in [x]^{\ell^p}_{\varepsilon'}
\]
\end{definition}

By definitions of $x \sim^{\ell^p}_{-\varepsilon +} y$ and $[x]^{\ell^p}_{-\varepsilon+}$, we have the following equivalent. 

\begin{lemma}
The following conditions are equivalent for any $p \in [1,\infty]$ and any number $\varepsilon \geq 0$ and any points $x,y \in X$: 
\\
{\rm(1)} $x \sim^{\ell^p}_{-\varepsilon +} y$.
\\
{\rm(2)}  $y \in \bigcap_{\varepsilon' \in [0,\varepsilon]} \bigcap_{z \in [x]^{\ell^p}_{\varepsilon'}}[z]^{\ell^p}_{\varepsilon'+}$. 
\\
{\rm(3)} $y \in [x]^{\ell^p}_{-\varepsilon+}$. 
\end{lemma}

We introduce the following concept, which is ``chain recurrence with negative errors''. 

\begin{definition}\label{def:cr02}
For any $p \in [1,\infty]$ and any non-negative number $\varepsilon \geq 0$, a point $x \in X$ is {\bf $\bm{\varepsilon}$-non-$\bm{\ell^p}$-gradient} if $x \in [x]^{\ell^p}_{-\varepsilon+}$. 
Denote by $\bm{\mathrm{CR}^{\ell^p}_{-\varepsilon}(f)}$ the set of $\varepsilon$-non-$\ell^p$-gradient points. 
\end{definition}

\begin{definition}\label{def:cr2_neg}
The subset $\mathrm{CR}^{\ell^\infty}_{-\varepsilon}(f)$ is denoted by $\bm{\mathrm{CR}_{-\varepsilon}(f)}$ and its points are called {\bf $\bm{-\varepsilon}$-chain recurrent points}.
Similarly, the subset $\mathrm{CR}^{\ell^1}_{-\varepsilon}(f)$ is denoted by $\bm{\mathrm{CR}^\Sigma_{-\varepsilon}(f)}$ and its points are called {\bf $\bm{-\varepsilon}$-$\bm{\sum}$chain recurrent points}.
\end{definition}

Roughly speaking, the previous definition formalizes the following concept: 
In a recurrent region, one can return to the starting point using approximately the same amount of energy as was used to move away from it. However, upon entering a gradient region, it becomes impossible to return using nearly the same amount of energy.
As above, we have the following observation. 

\begin{lemma}
The following conditions are equivalent for any $p \in [1,\infty]$ and any number $\varepsilon \geq 0$ and any point $x \in X$: 
\\
{\rm(1)} $x$ is $\varepsilon$-non-$\ell^p$-gradient. 
\\
{\rm(2)}  $x \in \bigcap_{\varepsilon' \in [0,\varepsilon]} \bigcap_{z \in [x]^{\ell^p}_{\varepsilon'}}[z]^{\ell^p}_{\varepsilon'+}$. 
\\
{\rm(3)} $x \sim^{\ell^p}_{-\varepsilon +} x$.
\end{lemma}

We define the $-\varepsilon$-chain recurrent components to define the Morse (hyper-)graphs with $\varepsilon$-$\ell^p$-errors (see Definitions~\ref{def:MG_error} and \ref{def:Morse_hyp_graph_map02} for details).  

\begin{definition}\label{def:cr04}
For any $p \in [1,\infty]$ and any non-negative number $\varepsilon \geq 0$, define a binary relation $\approx^{\ell^p}_{\mathrm{cr}_{-\varepsilon+}}$ on $\mathrm{CR}^{\ell^p}_{-\varepsilon}(f)$ by $x \approx^{\ell^p}_{\mathrm{cr}_{-\varepsilon+}} y$ if $x \sim^{\ell^p}_{-\varepsilon+} y$ and $y \sim^{\ell^p}_{-\varepsilon+} x$.
The transitive closure of $\approx_{\mathrm{cr}^{\ell^p}_{-\varepsilon +}}$ on $\mathrm{CR}^{\ell^p}_{-\varepsilon}(f)$ is an equivalence relation on $\mathrm{CR}^{\ell^p}_{-\varepsilon}(f)$, denote by $\approx_{\mathrm{CR}^{\ell^p}_{-\varepsilon +}}$. 
We call an equivalence class of the chain recurrent set $\mathrm{CR}^{\ell^p}_{-\varepsilon}(f)$ with respect to $\approx_{\mathrm{CR}^{\ell^p}_{-\varepsilon+}}$ a {\bf $\bm{-\varepsilon}$-chain recurrent component} of $\mathrm{CR}^{\ell^p}_{-\varepsilon}(f)$. 
\end{definition}

We obtain the following inclusions. 

\begin{lemma}\label{lem:41}
The following properties hold for any $p \in [1,\infty]$: 
\\
{\rm(1)} $\mathrm{CR}^{\ell^p}_{-0}(f) \subseteq \mathrm{CR}^{\ell^p}_{0}(f)$. 
\\
{\rm(2)} For any non-negative numbers $\varepsilon_1 < \varepsilon_2$, we have $\mathrm{CR}^{\ell^p}_{-\varepsilon_2}(f) \subseteq \mathrm{CR}^{\ell^p}_{-\varepsilon_1}(f)$. 
\\
{\rm(3)} If $\mathrm{CR}^{\ell^p}(f) = X$, then $\mathrm{CR}^{\ell^p}_{-\varepsilon}(f) = X$ for any non-negative number $\varepsilon \geq 0$. 
\end{lemma}

\begin{proof}
Recall that $\mathrm{CR}^{\ell^p}_{0}(f)$ is the set of $0$-$\ell^p$-chain recurrent points. 
In other words, we have $\mathrm{CR}^{\ell^p}_{0}(f) = \{ x \in X \mid x \sim^{\ell^p}_{0 +} x\}$. 
Therefore, we have $\mathrm{CR}^{\ell^p}_{-0}(f) = \{ x \in X \mid y \sim^{\ell^p}_{0 +} x \text{ for any }y \in O^+(x) \} \subseteq \{ x \in X \mid x \sim^{\ell^p}_{0 +} x\} = \mathrm{CR}^{\ell^p}_{0}(f)$. 
This means that assertion (1) holds. 
By definition of $\mathrm{CR}^{\ell^p}_{-\varepsilon}(f)$, assertions (2) and (3) hold. 
\end{proof}

\subsubsection{Total order on a variance of $\R$}

Consider the pair $(\R, \leq_{\R})$ of the real line and the standard order. 
Write two distinct points $-0, +0$ which are not contained in $\R$. 
Set $(-\infty,-0] := (\infty,0) \sqcup \{ -0 \}$ and $[0,\infty) := \{ +0 \} \sqcup (0,\infty)$. 
By considering $-0$ as the maximal element in $(-\infty,-0]$ and $+0$ as the minimal element in $[0,\infty)$, the subsets $(-\infty,-0]$ and $[0,\infty)$ are totally ordered. 
By setting $-0<+0$, the disjoint union $(-\infty,-0] \sqcup [0,\infty)$ becomes a totally ordered set. 

\subsubsection{Filtration whose parameters forms $(-\infty,-0] \sqcup [0,\infty)$}

Lemma~\ref{lem:equality} and the previous lemma imply the following inclusions. 

\begin{proposition}\label{th:main01}
For any $p \in [1,\infty]$, the family $(\mathrm{CR}^{\ell^p}_{\varepsilon}(f))_{\varepsilon \in (-\infty,-0] \sqcup [0,\infty)}$ of subsets of $X$ is a filtration of $X$. 
\end{proposition}

\subsubsection{Correspondence between $\mathrm{CR}^{\ell^p}_{0}$ and $\mathrm{CR}^{\ell^p}_{-0}$ under continuity}

To show the correspondence between $\mathrm{CR}^{\ell^p}_{0}$ and $\mathrm{CR}^{\ell^p}_{-0}$ under continuity, we have the following statement.

\begin{lemma}\label{lem:forward_inv}
Let $f \colon X \to X$ be a continuous mapping on a metric space $X$. 
For any $p \in [1,\infty]$ and any $x,x' \in X$ with $x' \sim^{\ell^p}_{0 +} x$, we have $f(x') \sim^{\ell^p}_{0 +} x$. 
\end{lemma}

\begin{proof}
Fix any positive number $\varepsilon$. 
If $d(O^+(f(x')), x) < \varepsilon$, then $f(x') \sim^{\ell^p}_{\varepsilon} x$. 
Thus we may assume that $d(O^+(f(x')), x) \geq \varepsilon$. 
Choose a positive number $\delta_2  \in (0, \varepsilon/2)$ with $f(B_{\delta_2}(f(x')) \subset B_{\varepsilon/4}(f^2(x'))$. 
Set a number $\delta_1 \in (0, \min\{ \delta_2/2, \varepsilon/4 \})$. 
By $x' \sim^{\ell^p}_{0 +} x$, there is a $\delta_1$-$\ell^p$-chain $(x', y_2, \ldots , y_{n-1}, y_n = x)$. 
From $d(O^+(f(x')), x) \geq \varepsilon > \delta_1$, we have $y_2 \neq x$ and so $n \geq 3$. 
Consider the sequence $(f(x'), y_3, \ldots , y_{n-1}, x)$. 
From $d(y_2,f(x'))< \delta_1 < \delta_2$, we obtain $y_2 \in B_{\delta_2}(f(x'))$. 
By $f(B_{\delta_2}(f(x')) \subset B_{\varepsilon/4}(f^2(x'))$, we have $d(f^2(x'), y_3) \leq d(f^2(x'), f(y_2)) + d(f(y_2), y_3) \leq \varepsilon/4 + \delta_1 < \varepsilon/2$. 
For any $p \in [1, \infty)$, since 
\[
\left( \left( \dfrac{\varepsilon}{2} \right)^p + \delta_1^p \right)^{\dfrac{1}{p}} < \left( \left( \dfrac{\varepsilon}{2} \right)^p + \left( \dfrac{\varepsilon}{4} \right)^p \right)^{\dfrac{1}{p}} < \left( 2 \left( \dfrac{\varepsilon}{2} \right)^p \right)^{\dfrac{1}{p}} \leq 2 \dfrac{\varepsilon}{2} = \varepsilon,
\]
the sequence $(f(x'), y_3, \ldots , y_{n-1}, x)$ is an $\varepsilon$-$\ell^p$-chain. 
If $p = \infty$, then the inequality $\max \{ \delta_1, \varepsilon/2 \} = \varepsilon/2 < \varepsilon$ implies that 
$(f(x'), y_3, \ldots , y_{n-1}, x)$ is an $\varepsilon$-$\ell^\infty$-chain. 
\end{proof}

We have the correspondence. 

\begin{lemma}\label{lem:-0=0}
Suppose that $f \colon X \to X$ is a continuous mapping. 
For any $p \in [1,\infty]$, we have that $\mathrm{CR}^{\ell^p}_{-0}(f) = \mathrm{CR}^{\ell^p}_{0}(f) = \mathrm{CR}^{\ell^p}(f)$.
\end{lemma}

\begin{proof}
Lemma~\ref{lem:41} implies $\mathrm{CR}^{\ell^p}_{-0}(f) \subseteq \mathrm{CR}^{\ell^p}_{0}(f)$. 
Notice that $\mathrm{CR}^{\ell^p}_{0}(f) = \mathrm{CR}^{\ell^p}(f)$. 
Fix a point $x \in \mathrm{CR}^{\ell^p}_{0}(f) = \mathrm{CR}^{\ell^p}(f)$.
Then $x \sim^{\ell^p}_{0 +} x$. 
Fix a point $y \in [x]^{\ell^p}_{0} = \{ z \in X \mid x \sim^{\ell^p}_0 z \} = O^+(x)$.
Then there is a positive integer $k \in \Z_{>0}$ such that $f^k(x) = y$. 
By applying Lemma~\ref{lem:forward_inv} $k$ times, we obtain that $y = f^k(x) \sim^{\ell^p}_{0+} x$. 
This means that $x \in \mathrm{CR}^{\ell^p}_{-0}(f)$. 
\end{proof}

Identifying $-0$ and $+0$, one can identify the resulting totally ordered set of $(-\infty,-0] \sqcup [0,\infty)$ with $\R$. 
Using this identification, Proposition~\ref{th:main01} and Lemma~\ref{lem:-0=0} imply the following statement. 

\begin{theorem}\label{thm:01_gen}
For any $p \in [1, \infty]$ and any continuous mapping $f \colon X \to X$ on a metric space $X$, the family $(\mathrm{CR}^{\ell^p}_{\varepsilon}(f))_{\varepsilon \in \R}$ is a filtration of $X$ with $\mathrm{CR}^{\ell^p}(f) = \mathrm{CR}^{\ell^p}_{0}(f)$. 
\end{theorem}

The previous result implies Theorem~\ref{main:01} and Theorem~\ref{main:02}. 
On the other hand, we have the following non-correspondence. 

\begin{lemma}
For any $p \in [1, \infty]$, there is a non-continuous bijection $f \colon X \to X$ with $\mathrm{CR}^{\ell^p}_{-0}(f) \subsetneq \mathrm{CR}^{\ell^p}_{0}(f)$.
\end{lemma}

\begin{proof}
Put $A := \{0 \} \times (0,1]$, $B := \{1 \} \times (0,1]$, $C := \Z \times \{ 0 \}$, and $X := A \sqcup B \sqcup C \subset \R^2$. 
Define a bijection $f \colon X \to X$ by $f \vert_{C}(x,y) = (x+1,y)$, $f \vert_{A}(0,y) = (1,y)$, and $f \vert_{B}(1,y) = (0,y)$.
Then $(0,0) \sim^{\ell^p}_0 (2,0)$ and $(2,0) \not\sim^{\ell^p}_{0+} (0,0)$. 
However, we have $(0,0) \sim^{\ell^p}_0 (1,0)\sim^{\ell^p}_{0+} (0,0)$. 
Therefore, we obtain $x \in \mathrm{CR}^{\ell^p}_{0}(f) - \mathrm{CR}^{\ell^p}_{-0}(f)$. 
\end{proof}

\subsection{Singular limits} 

We have the following singular limit behavior of $\mathrm{CR}_{-\varepsilon}^\Sigma(f)$.

\begin{theorem}\label{th:counter_ex04}
For any $\varepsilon>0$, there is a homeomorphism $f$ on a circle with $\bigcup_{\varepsilon>0} \mathrm{CR}_{-\varepsilon}^\Sigma(f) \subsetneq \mathrm{CR}_{-0}^\Sigma(f) = \mathrm{CR}_{0}^\Sigma(f)$. 
\end{theorem}

In fact, there is such a homeomorphism in \S~\ref{sec:c_ex_001}. 
\begin{figure}[t]
\begin{center}
\includegraphics[scale=0.345]{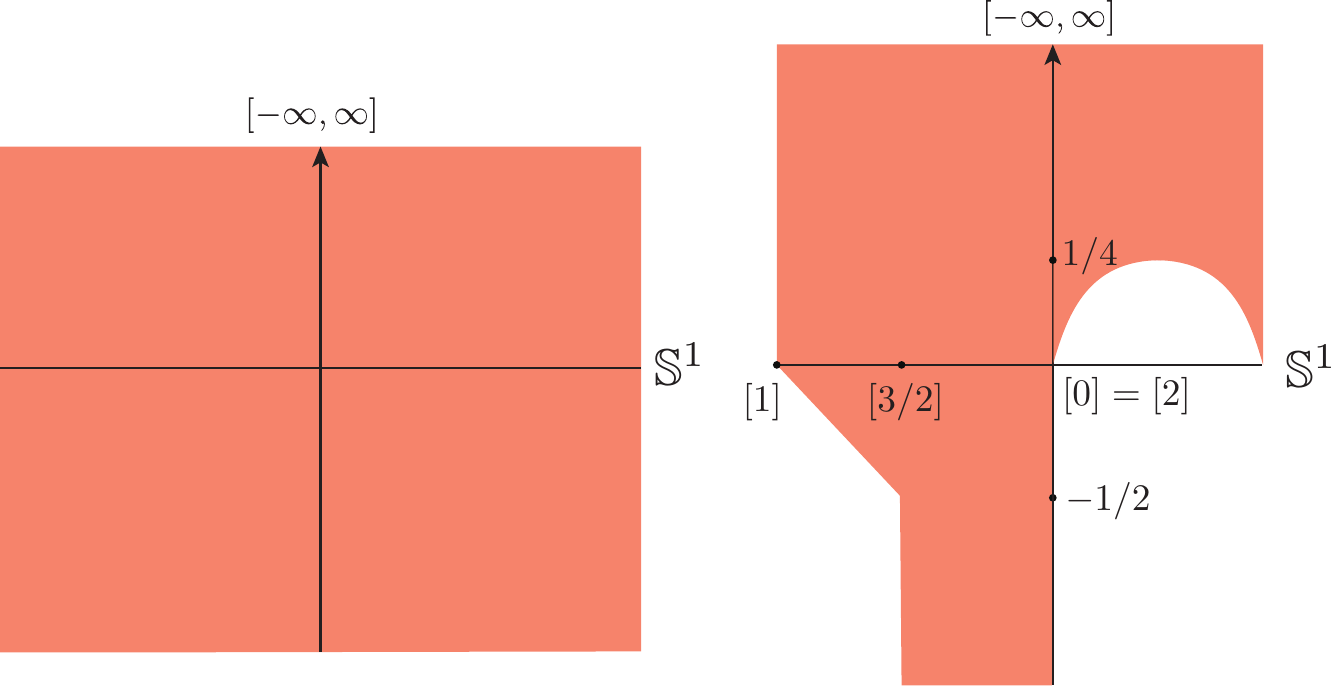}
\end{center} 
\caption{Diagrams $D_{\mathrm{CR}}(f_{R})$ and $D_{\mathrm{CR}}^\Sigma(f_{R})$.}
\label{fig:ex_diagram_00}
\end{figure} 
Notice that the subsets $\mathrm{CR}_{-\varepsilon}(f)$ and $\mathrm{CR}_{-\varepsilon}^\Sigma(f)$ are incomparable in general. 
In fact, a homeomorphism $f_{\mathrm{rep}} \colon \R \to \R$ defined by 
\[
f_{\mathrm{rep}}(x) = 
\begin{cases}
x &  (x \leq 0) \\
2x &  (x > 0) 
\end{cases}
\]
satisfies that $\mathrm{CR}_{-\varepsilon}(f) = \emptyset \subsetneq \R_{\leq - \varepsilon} = \mathrm{CR}_{-\varepsilon}^\Sigma(f)$ for any $\varepsilon > 0$.  
Moreover, a homeomorphism $f_R \colon \R/2\Z \to \R/2\Z$ defined by 
\[
f_{R}([x]) = 
\begin{cases}
[x^2] &  (x \in [0,1]) \\
[x] &  (x \in [1,2]) 
\end{cases}
\]
satisfies that $\mathrm{CR}_{\varepsilon}(f_{R}) = \R/2\Z$ for any $\varepsilon \in \R$ and 
\[
\mathrm{CR}_{\varepsilon}^\Sigma(f_{R}) = 
\begin{cases}
\R/2\Z &  (\varepsilon >1/4) \\
\left[\dfrac{1+ \sqrt{1-4\varepsilon}}{2} ,2+ \dfrac{1- \sqrt{1-4\varepsilon}}{2} \right]/2\Z & (\varepsilon \in (0,1/4]) \\
[1 - \varepsilon,2]/2\Z &  (\varepsilon \in [-1/2,0]) \\
[3/2,2]/2\Z &  (\varepsilon < -1/2) 
\end{cases}
\]
and so that $\mathrm{CR}_{-\varepsilon}^\Sigma(f_{R}) \subsetneq \mathrm{CR}_{-\varepsilon}(f_{R})$ for any $\varepsilon>0$. 

\subsubsection{$\varepsilon$-$\ell^p$-chain recurrent potential and $\varepsilon$-$\ell^p$-chain recurrent diagram of a mapping}

As bifurcation diagrams of dynamical systems, we define the $\ell^p$-chain recurrent potential and the $\varepsilon$-$\ell^p$-chain recurrent diagram to analyze dynamical systems as follows. 

Let $f \colon X \to X$ be a mapping (resp. continuous mapping) on a metric space $X$ and put $I := [-\infty,-0] \sqcup [0,\infty]$ (resp. $I := [-\infty,\infty]$). 
Write $\mathrm{CR}^{\ell^p}_{\infty}(f) := X$ for any $p \in [1,\infty]$.

\begin{definition}\label{def:debut}
For any $p \in [1,\infty]$, define the {\bf $\ell^p$-chain recurrent potentials} \(\tau^{\ell^p}_f \colon X \rightarrow I\) as follows:
\[
\tau^{\ell^p}_f(x) := \inf \{\varepsilon \in I - \{-\infty, \infty \} \mid x \in \mathrm{CR}^{\ell^p}_{\varepsilon}(f) \}
\]
\end{definition}

\begin{definition}\label{def:crd_map_001}
For any $p \in [1,\infty]$, define the {\bf $\ell^p$-chain recurrent diagram} $D^{\ell^p}_{\mathrm{CR}}(f)$ of $f$ as follows: 
\\
\[
D^{\ell^p}_{\mathrm{CR}}(f) := \bigsqcup_{\varepsilon \in I} \mathrm{CR}^{\ell^p}_{\varepsilon}(f) \times \{ \varepsilon \} = \{(x,\varepsilon) \in X \times I \mid \tau^{\ell^p}_f(x) \leq \varepsilon \}
\]
\end{definition}

\begin{definition}\label{def:crd_map_002}
The function $\tau^{\ell^\infty}_f$ is denoted by $\bm{\tau_f}$ and is called the {\bf chain recurrent potential}.
The subset $D^{\ell^\infty}_{\mathrm{CR}}(f)$ is denoted by $\bm{D_{\mathrm{CR}}(f)}$ and is called the  {\bf chain recurrent diagram}. 

Similarly, the function $\tau^{\ell^1}_f$ is denoted by $\bm{\tau^\Sigma_f}$ and is called the {\bf $\bm{\sum}$chain recurrent potential}.
The subset $D^{\ell^1}_{\mathrm{CR}}(f)$ is denoted by $\bm{D^\Sigma_{\mathrm{CR}}(f)}$ and is called the  {\bf $\bm{\sum}$chain recurrent diagram}. 
\end{definition}

For instance, the chain recurrent potential, the $\sum$chain recurrent potential, the $\sum$chain recurrent diagram, and chain recurrent diagram of the mapping $f_R$ are 
\[
\tau_{f_R} = -\infty
\]
\[
\tau^\Sigma_{f_R}([x]) = \begin{cases}
x- x^2 & (x \in (0,1)) \\
1 - x &  (x \in [1,3/2)) \\
- \infty &  (x \in [3/2,2]) 
\end{cases}
\]
\[
D_{\mathrm{CR}}(f_R) = (\R/2\Z) \times [-\infty, \infty]
\]
\[
\begin{split}
D_{\mathrm{CR}}^\Sigma(f_R) &= (([3/2,2]/2\Z) \times [-\infty, -1/2]) 
\\
& \hspace{8pt}  \sqcup \{([x],\varepsilon) \mid x \in [1 - \varepsilon,2], \varepsilon \in (-1/2,0] \}
\\
& \hspace{8pt} \sqcup \left\{([x],\varepsilon)  \middle| x \in \left[\dfrac{1+ \sqrt{1-4\varepsilon}}{2} ,2+ \dfrac{1- \sqrt{1-4\varepsilon}}{2} \right], \varepsilon \in (0, 1/4] \right\}
\\
& \hspace{8pt} \sqcup (\R/2\Z \times (1/4, \infty])
\end{split}
\]
as shown in Figure~\ref{fig:ex_diagram_00}. 

Moreover, consider a contraction $g \colon \R \to \R$ by $g(x)= x/2$. 
Then the chain recurrent potential and the $\ell^p$-chain recurrent potential \(\tau^{\ell^p}_g \colon \R \rightarrow I\) are
\[
\tau^{\ell^p}_g(x) = 
\begin{cases}
-\infty &  (x = 0) \\
\vert x/2 \vert &  (x \neq 0) 
\end{cases}
\]
and the $\ell^p$-chain recurrent diagram of $g$ are 
\[
D^{\ell^p}_{\mathrm{CR}}(g) = (\{ 0 \} \times [-\infty, 0)) \sqcup \{ (x, \varepsilon) \mid \varepsilon \geq 0, x \in [-\varepsilon/2, \varepsilon/2]\}
\]
and the $\ell^p$-circulation cost of $g$ are $\infty$. 
\begin{figure}[t]
\begin{center}
\includegraphics[scale=0.345]{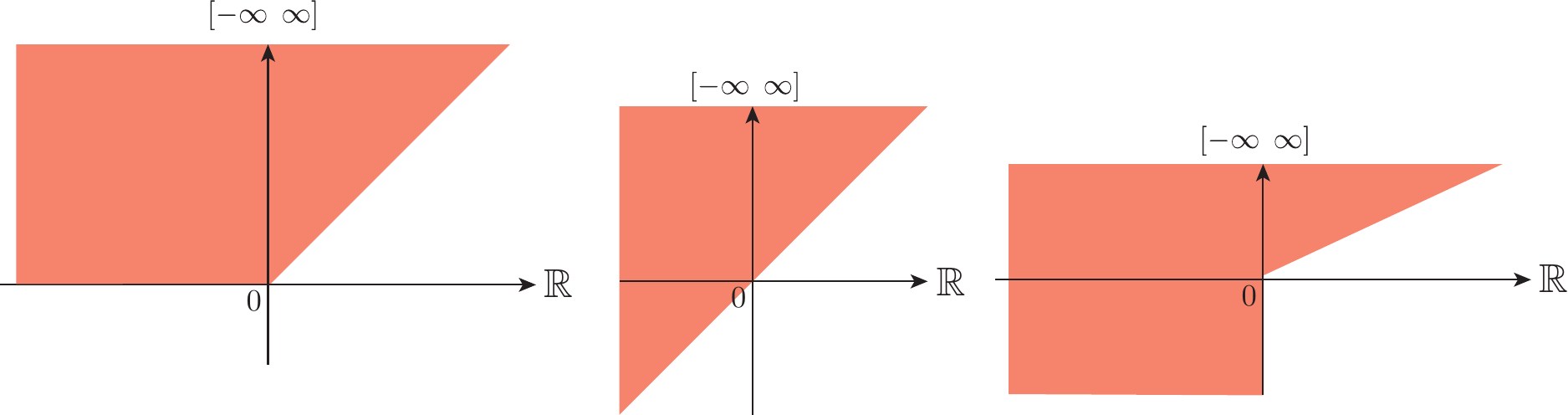}
\end{center} 
\caption{Diagrams $D_{\mathrm{CR}}(f_{\mathrm{rep}})$, $D_{\mathrm{CR}}^\Sigma(f_{\mathrm{rep}})$, and $D^{\ell^p}_{\mathrm{CR}}(f_{\mathrm{att}})$.}
\label{fig:ex_diagram}
\end{figure} 

Furthermore, consider the homeomorphism $f_{\mathrm{rep}}$ as above, and a homeomorphism $f_{\mathrm{att}} \colon \R \to \R$ defined as follows: 
\[
f_{\mathrm{att}}(x) = 
\begin{cases}
x &  (x \leq 0) \\
x/2 &  (x > 0) 
\end{cases}
\]
Then the chain recurrent potential $\tau_{f_{\mathrm{rep}}}$, the $\sum$chain recurrent potential $\tau_{f_{\mathrm{rep}}}^{\Sigma}$, and the $\ell^p$-chain recurrent potential \(\tau^{\ell^p}_{f_{\mathrm{att}}} \colon \R \rightarrow [-\infty, \infty) \) are
\[
\tau_{f_{\mathrm{rep}}}(x) = 
\begin{cases}
0 &  (x \leq 0) \\
x &  (x > 0) 
\end{cases}
\]
\[
\tau_{f_{\mathrm{rep}}}^{\Sigma}(x) = x
\]
\[
\tau^{\ell^p}_{f_{\mathrm{att}}}(x) = 
\begin{cases}
-\infty &  (x \leq 0) \\
x/2 &  (x > 0) 
\end{cases}
\]
and the chain recurrent diagram and the $\sum$chain recurrent diagram of $f_{\mathrm{rep}}$ and $\ell^p$-chain recurrent diagram of  $f_{\mathrm{att}}$ are 
\[
D_{\mathrm{CR}}(f_{\mathrm{rep}}) = \{ (x, \varepsilon) \mid \max \{0, x \} \leq \varepsilon\}
\]
\[
D_{\mathrm{CR}}^\Sigma(f_{\mathrm{rep}}) = \{ (x, \varepsilon) \mid x \leq \varepsilon \}
\]
\[
D_{\mathrm{CR}}^{\ell^p}(f_{\mathrm{att}}) = (\R_{\leq 0} \times [-\infty,\infty]) \sqcup \{ (x, \varepsilon) \mid 0 < x/2 \leq \varepsilon \}
\]
and the circulation $\sum$cost and the circulation cost of $f_{\mathrm{rep}}$ and the $\ell^p$-circulation cost $f_{\mathrm{att}}$ are $\infty$ as shown in Figure~\ref{fig:ex_diagram}.
The diagrams of these homeomorphisms suggest that the negative part of the filtration $D^{\ell^p}_{\mathrm{CR}}$ illustrates a part of the behavior near the $\ell^p$-chain recurrent set.

\begin{figure}[t]
\begin{center}
\includegraphics[scale=0.5]{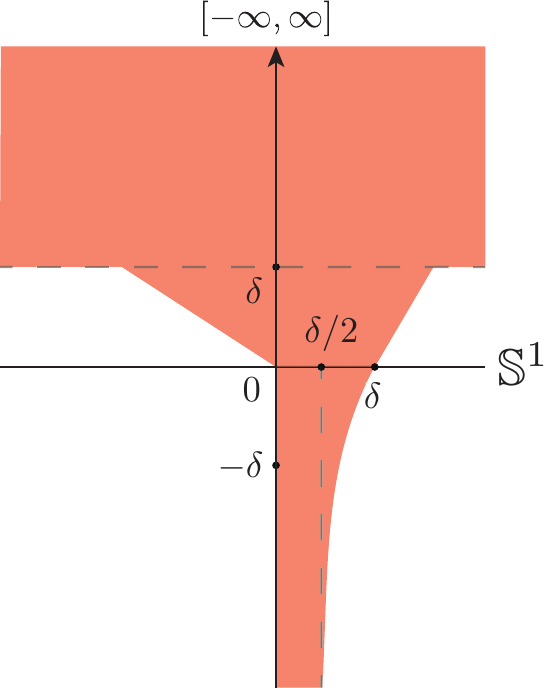}
\end{center} 
\caption{A schematic picture of the $\sum$chain recurrent diagram of a mapping $f$ defined in \S~\ref{sec:c_ex_001}.}
\label{fig:csr_diagram_01}
\end{figure}

In addition, Figure~\ref{fig:csr_diagram_01} is an example of a schematic picture of the $\sum$chain recurrent diagram of a mapping $f$ defined in \S~\ref{sec:c_ex_001}.

\section{Morse (hyper-)graph with finite errors for homeomorphisms}

In this section, we introduce Morse (hyper-)graphs with finite $\ell^p$-errors for homeomorphisms.  
%

\subsection{Morse graphs with $\varepsilon$-$\ell^p$-errors for homeomorphisms}\label{MG_homeo}

We introduce the Morse graphs with $\varepsilon$-$\ell^p$-errors for homeomorphisms as follows.

\begin{definition}\label{def:MG_error}
For any $p \in [1,\infty]$ and any homeomorphism $f$ on a metric space, the Morse graph $G_{\mathcal{X}}(f)$ of $\mathcal{X}$ with persistent directed edges in \S~\ref{sec:Morse_graph} is the {\bf Morse (directed) graph} of $f$ with $\varepsilon$-$\ell^p$-errors  if $\mathcal{X}$ is the set of $\varepsilon$-chain recurrent component of $\mathrm{CR}^{\ell^p}_{\varepsilon}(f)$.
Then we denote by $\bm{G^{\ell^p}_{\varepsilon}(f)}$ the Morse graph of $f$ with $\varepsilon$-$\ell^p$-errors. 
\end{definition}

\begin{definition}\label{def:MG_error_02}
The graph $G_{\varepsilon}(f) := G^{\ell^\infty}_{\varepsilon}(f)$ is called the Morse graph of $f$ with $\varepsilon$-errors. 
Similarly, the graph $G^{\Sigma}_{\varepsilon}(f) := G^{\ell^1}_{\varepsilon}(f)$ is called the Morse graph of $f$ with $\varepsilon$-$\sum$errors. 
\end{definition}

Notice that $G^{\ell^p}_{\varepsilon}(f)$ is a generalization of the Morse graph $G(f)$ of $f$ such that $G^{\ell^\infty}_{0}(f) = G_{0}(f) = G(f)$. 
Note that we can define $G'^{\ell^p}_{\varepsilon}(f)$ by replacing $G_{\mathcal{X}}(f)$ with $G'_{\mathcal{X}}(f)$, though directed edges of $G'(f)$ need not be persistent as shown in \S~\ref{sec:non-increasing}. 

\subsection{Vertex collapse}

To state a one-parameter family of Morse graphs, we introduce a vertex collapse, which is a concept analogous to edge collapse.
Recall an edge collapsing operation as follows.

\begin{definition}
A directed graph $G =(V,D)$ is obtained from a directed graph $G' = (V',D')$ by an {\bf edge collapse} of a directed edge $(a_1',a_2') \in D'$ into a vertex $a \in V$ if $V - \{ a \} = V' - \{ a_1',a_2' \}$ and a mapping $f_D \colon D' - \{ (a_1',a_2') \} \to D$ defined by 
\[ f_D((x,y)) := 
\begin{cases}
(x,y) & (\{ x, y \} \cap \{ a_1',a_2' \} = \emptyset) \\
(a,y) & (x \in \{ a_1',a_2' \} \text{ and } y \notin \{ a_1',a_2' \}) \\
(a,a) & ((x, y) =  (a_2',a_1')) \\
(x,a) & ( x \notin \{ a_1',a_2' \} \text{ and } y \in \{ a_1',a_2' \}) 
\end{cases}
\]
is well-defined and surjective. 
\end{definition}

we define the following collapsing operation using transfinite induction. 

\begin{definition}\label{def:ec}
Let $G=(V,D)$ be a directed graph, $G'=(V',D')$ a directed graph, $f_V \colon V' \to V$ a mapping,  
 and $D'_{\mathrm{loop}} := \{ (x',y') \in D' \mid x' \neq y', f_V(x') = f_V(y') \}$. 
A directed graph $G$ is obtained (resp. partially obtained) from a directed graph $G'$ by {\bf vertex collapses} associated to the mapping $f_V$ if 
a mapping $f_D \colon D' - D'_{\mathrm{loop}} \to D$ defined by $f_D((x,y)) := (f_V(x),f_V(y))$ is well-defined and surjective (resp. well-defined). 
Then we denote by $G' \mathrel{\raisebox{.85ex}{\rotatebox{-45}{$\twoheadrightarrow$}}}^{\hspace{-3pt}f_V} G$ (resp. $G' \searrow^{\hspace{-3pt}f_V} G$). 
\end{definition}

Notice that the restriction $f_D \vert_{(D' - D'_{\mathrm{loop}}) \cap V'^2_{\mathrm{inj}}}$ is injective for any directed graph $G$ with $G' \searrow^{\hspace{-3pt}f_V} G$. 

\subsection{Graph collapses}

We have the following reduction. 

\begin{theorem}\label{th:graph_collapse_homeomorphism}
Let $f \colon X \to X$ be a homeomorphism on a metric space $X$. 
The following statements hold for any $p \in [1,\infty]$ and any $\varepsilon_1 < \varepsilon_2 \in (-\infty,-0] \sqcup [0,\infty)$:
\\
{\rm(1)} The Morse graph $G^{\ell^p}_{\varepsilon_2}(f) = (V_2, D_2)$ is partially obtained from $G^{\ell^p}_{\varepsilon_1}(f) = (V_1, D_1)$ by vertex collapses associated to $h_{\varepsilon_2}^{\varepsilon_1}$ {\rm(i.e.} $G^{\ell^p}_{\varepsilon_1}(f) \searrow^{\hspace{-3pt}h_{\varepsilon_2}^{\varepsilon_1}} G^{\ell^p}_{\varepsilon_2}(f)${\rm)}, where the mapping $h_{\varepsilon_2}^{\varepsilon_1} \colon V_1 \to V_2$ is defined by assigning to each $X_1 \in V_1$ the $\varepsilon_2$-chain recurrent component of $\mathrm{CR}^{\ell^p}_{\varepsilon_2}(f)$ containing $X_1$. 
\\
{\rm(2)} If $\varepsilon_1 \geq 0$, then $h_{\varepsilon_2,D}^{\varepsilon_1}$ is surjective {\rm(i.e.} $G^{\ell^p}_{\varepsilon_1}(f) \mathrel{\raisebox{.85ex}{\rotatebox{-45}{$\twoheadrightarrow$}}}^{h_{\varepsilon_2}^{\varepsilon_1}} G^{\ell^p}_{\varepsilon_2}(f)${\rm)}.
\end{theorem}

\begin{proof}
Fix any $\varepsilon_1 < \varepsilon_2 \in (-\infty,-0] \sqcup [0,\infty)$. 
By definition, we have that $\mathrm{CR}^{\ell^p}_{\varepsilon_1}(f) \subseteq \mathrm{CR}^{\ell^p}_{\varepsilon_2}(f)$. 
Therefore, for any $\varepsilon_1$-chain recurrent component $X_1$ of $\mathrm{CR}^{\ell^p}_{\varepsilon_1}(f)$, there are an element $\varepsilon_2 > \varepsilon_1$ and an $\varepsilon_2$-chain recurrent component $X_2$ of $\mathrm{CR}^{\ell^p}_{\varepsilon_2}(f)$ with $X_1 \subseteq X_2$, because of $\mathrm{CR}^{\ell^p}_{\varepsilon_1}(f) \subseteq \mathrm{CR}^{\ell^p}_{\varepsilon_2}(f)$. 
Therefore, the mapping $h_{\varepsilon_2}^{\varepsilon_1} \colon V_1 \to V_2$ is well-defined.
Apply vertex collapses associated to $h_{\varepsilon_2}^{\varepsilon_1}$, the resulting graph from $G^{\ell^p}_{\varepsilon_1}(f)$ is a spanning subgraph of $G^{\ell^p}_{\varepsilon_2}(f)$. 
By definition of persistent directed edgeset of $G^{\ell^p}_{\varepsilon}(f)$, since $\mathrm{CR}^{\ell^p}_{\varepsilon_1}(f) \subseteq \mathrm{CR}^{\ell^p}_{\varepsilon_2}(f)$, the map $h_{\varepsilon_2,D}^{\varepsilon_1} \colon D_1 \to D_2$ defined by 
\[
h_{\varepsilon_2,D}^{\varepsilon_1}(X_1,X'_1) := (h_{\varepsilon_2}^{\varepsilon_1}(X_1),h_{\varepsilon_2}^{\varepsilon_1}(X'_1))
\]
 is well-defined.
%
%

Suppose that $\varepsilon_1 \geq 0$. 
We have the following relations: 
\[
\bigcup_{x \in X} \omega(x) \cup \alpha(x) \subseteq \mathrm{CR}^{\ell^p}(f) \subseteq \mathrm{CR}^{\ell^p}_{\varepsilon_1}(f) \subseteq \mathrm{CR}^{\ell^p}_{\varepsilon_2}(f)
\]
Fix any $(X_j,X_k) \in D_2$. 
By definition of $D_2$, there are points $x, \alpha, \omega \in X$ such that $\alpha \in \alpha(x) \cap X_j$ and $\omega \in \omega(x) \cap X_k$. 
From $\alpha, \omega \in \bigcup_{x \in X} \omega(x) \cup \alpha(x) \subseteq \mathrm{CR}^{\ell^p}_{\varepsilon_1}(f)$, there are vertices $X'_{j'}, X'_{k'} \in V_1$ such that $\alpha \in X'_{j'} \subseteq X_j$ and $\omega \in X'_{k'} \subseteq X_k$. 
Then $\alpha \in \alpha(x) \cap X'_{j'}$ and $\omega \in \omega(x) \cap X'_{k'}$. 
This means that $(X'_{j'}, X'_{k'}) \in D_1$. 
By definition of $h_{\varepsilon_2,D}^{\varepsilon_1}$, we have $h((X'_{j'}, X'_{k'})) = (X_j,X_k)$. 
Therefore, the map $h_{\varepsilon_2,D}^{\varepsilon_1} \colon D_1 \to D_2$ is surjective.
\end{proof}

Notice that the family $(h_{\varepsilon}^{\varepsilon'})_{\varepsilon' < \varepsilon}$ as in the previous theorem satisfies the cocycle condition (i.e. $h_{\varepsilon_3}^{\varepsilon_2} \circ h_{\varepsilon_2}^{\varepsilon_1} = h_{\varepsilon_3}^{\varepsilon_1}$ for any $\varepsilon_1 < \varepsilon_2 < \varepsilon_3 \in (-\infty,-0] \sqcup [0,\infty)$). 
%
%

\subsection{``Connecting'' properties and their coarse variants}

To define Morse graphs with finite errors for mappings, we weaken ``connecting'' properties as follows. 
We have the following properties. 

\begin{lemma}\label{lem:alpha}
Suppose that $f$ is a homeomorphism on a metric space $X$. 
The following statements hold for any $p \in [1,\infty]$, any point \( x \in X \), and any subset \( A \subseteq X \):
\\
{\rm(1)} $x \in \bigcup_{\alpha \in A} [\alpha]_{0+}^\Sigma$ implies $x \in \bigcup_{\alpha \in A} [\alpha]^{\ell^p}_{0+}
$.  
\\
{\rm(2)} 
$\alpha(x) \cap A \neq \emptyset$ implies $x \in \bigcup_{\alpha \in A} [\alpha]_{0+}^\Sigma$. 
\end{lemma}

\begin{proof}
For any $p \in [1,\infty)$ and any non-negative numbers $\varepsilon_0, \varepsilon_1, \ldots , \varepsilon_{n-1}$, since 
\[
(\sum_{i=0}^{n-1} \varepsilon_i)^p \geq \sum_{i=0}^{n-1} \varepsilon_i^p = (\| ( \varepsilon_i )_{i=0}^{n-1} \|_{p})^p
\]
 and 
 \[
 \| ( \varepsilon_i )_{i=0}^{n-1} \|_{\infty} = \max \{\varepsilon_0, \varepsilon_1, \ldots , \varepsilon_{n-1}  \},\]
 we have the following inequality: 
\[
\| ( \varepsilon_i )_{i=0}^{n-1} \|_{1} = \sum_{i=0}^{n-1} \varepsilon_i \geq \max \left\{  \| ( \varepsilon_i )_{i=0}^{n-1} \|_{p}, \| ( \varepsilon_i )_{i=0}^{n-1} \|_{\infty} \right\}
\]
This implies the assertion (1) for any $p \in [1,\infty]$. 

Suppose that $\alpha(x) \cap A \neq \emptyset$. 
If \(O^-(x) \cap A \neq \emptyset\), then there is a point $a' \in O^-(x) \cap A$ and so $x \in O^+(a') = [a']_0^\Sigma \subseteq [a']_{0+}^\Sigma \subseteq \bigcup_{\alpha \in A} [\alpha]_{0+}^\Sigma$.
Thus we may assume that \( O^-(x) \cap A = \emptyset\). 
Fix a point $\alpha \in \alpha(x) \cap A$. 
By $O^-(x) \cap A = \emptyset$, we obtain $\alpha \in (\alpha(x) \setminus O^-(x)) \cap A$. 
Then $U \cap O^-(x) \neq \emptyset$ for any open \nbd $U$ of $\alpha$. 
In particular, we have $B_\varepsilon(\alpha) \cap O^-(x) \neq \emptyset$ for any $\varepsilon > 0$. 
The continuity of $f$ implies that $\alpha \sim_{\varepsilon_\Sigma} x$ for any $\varepsilon > 0$.
Thus $\alpha \sim_{0_\Sigma +} x$ and so $x \in [\alpha]_{0+}^\Sigma \subseteq \bigcup_{\alpha \in A} [\alpha]_{0+}^\Sigma$.
%
\end{proof}

\begin{lemma}\label{lem:omega}
Let $f$ be a continuous mapping on a metric space $X$. 
The following statements hold for any $p \in [1,\infty]$, any point \( x \in X \), and any subset \( W \subseteq X \):
\\
{\rm(1)} $[x]_{0+}^\Sigma \cap W \neq \emptyset$ implies $[x]^{\ell^p}_{0+} \cap W \neq \emptyset$.
\\
{\rm(2)} 
$\omega(x) \cap W \neq \emptyset$ implies $[x]_{0+}^\Sigma \cap W \neq \emptyset$. 
\end{lemma}

\begin{proof}
As the proof of the previous lemma, for any $p \in [1,\infty)$ and any non-negative numbers $\varepsilon_0, \varepsilon_1, \ldots , \varepsilon_{n-1}$, we have the following inequality: 
\[
\sum_{i=0}^{n-1} \varepsilon_i \geq \max \left\{  \| ( \varepsilon_i )_{i=0}^{n-1} \|_{p}, \| ( \varepsilon_i )_{i=0}^{n-1} \|_{\infty} \right\}
\]
Therefore, the condition $[x]_{0+}^\Sigma \cap W \neq \emptyset$ implies $[x]^{\ell^p}_{0+} \cap W \neq \emptyset$  for any $p \in [1,\infty]$. 

Suppose that $\omega(x) \cap W \neq \emptyset$.
If \(O^+(x) \cap W \neq \emptyset\), then $\emptyset \neq O^+(x) \cap W = [x]_{0}^\Sigma \cap W \subseteq [x]_{0+}^\Sigma \cap W$. 
Thus we may assume that \(O^+(x) \cap W = \emptyset\). 
Fix a point $\omega \in \omega(x) \cap W$. 
By $O(x)^+ \cap W = \emptyset$, we obtain $\omega \in (\omega(x) \setminus O^+(x)) \cap W$.
Then $U \cap O^{+}(x) \neq \emptyset$ for any open \nbd $U$ of $\omega$.
In particular, we have $B_\varepsilon(\omega) \cap O^{+}(x) \neq \emptyset$ for any $\varepsilon > 0$. 
The continuity of $f$ implies that $x \sim_{\varepsilon_\Sigma} \omega$ for any $\varepsilon > 0$.
Thus $x \sim_{0_\Sigma +} \omega$ and so $\omega \in [x]_{0+}^\Sigma$. 
%
\end{proof}

\section{Morse (hyper-)graphs with $(\varepsilon,\nu)$-errors}

We introduce a coarse generalization applicable to arbitrary mappings, which naturally admits a canonical filtration.
Let $f \colon X \to X$ be a mapping on a metric space $X$. 
Fix a set $\mathcal{X} = \{ X_i \}_{i \in \Lambda}$ of disjoint subsets. \
We introduce Morse graphs with errors as follows. 

\begin{definition}\label{def:Morse_graph_map_001}
For any $p \in [1,\infty]$ and any number $\nu \geq 0$, a directed graph $(V, D(\nu))$ with the vertex set $V := \{ X_i \mid i  \in \Lambda \}$, and with the {\bf directed edge set} $D(\nu) := \{ (X_j, X_k) \mid D_{j,k}(\nu) \neq \emptyset, j \neq k \in \Lambda \}$ is called the {\bf Morse graph} of $\mathcal{X}$ with $\nu$-errors, where
\[
D_{j,k}(\nu) := \left\{ x \in X \middle| \, x \in \bigcup_{\alpha \in X_j} [\alpha]^{\ell^p}_{\nu+}, \,[x]^{\ell^p}_{\nu+} \cap X_k \neq \emptyset \right\}.
\]
Then such a graph is denoted by $\bm{G^{\ell^p}_{\mathcal{X},\nu}}$. 
\end{definition}

\begin{definition}\label{def:Morse_graph_map02}
For any $p \in [1,\infty]$ and any $\varepsilon \in (-\infty,-0] \sqcup [0,\infty)$ and any number $\nu \geq 0$, the graph $G^{\ell^p}_{\mathcal{X},\nu}$ is the {\bf Morse graph} of $f$ with $\ell^p$-$(\varepsilon,\nu)$-errors if $\mathcal{X}$ is the set of $\varepsilon$-chain recurrent components of $\mathrm{CR}^{\ell^p}_{\varepsilon}(f)$. 
Then we denote by $\bm{G^{\ell^p}_{(\varepsilon,\nu)}(f)}$ the Morse graph of $f$ with $(\varepsilon,\nu)$-$\ell^p$-errors. 
\end{definition}

\begin{definition}\label{def:MG_error_03}
The graph $\bm{G_{(\varepsilon,\nu)}(f)} := G^{\ell^\infty}_{(\varepsilon,\nu)}(f)$ is called the Morse graph of $f$ with $(\varepsilon,\nu)$-errors. 
Similarly, the graph $\bm{G^{\Sigma}_{(\varepsilon,\nu)}(f)} := G^{\ell^1}_{(\varepsilon,\nu)}(f)$ is called the Morse graph of $f$ with $(\varepsilon,\nu)$-$\sum$errors. 
\end{definition}

Lemma~\ref{lem:alpha} and Lemma~\ref{lem:omega} imply the following observation. 

\begin{lemma}
Let $f$ be a homeomorphism on a metric space. 
For any $\varepsilon \in (-\infty,-0] \sqcup [0,\infty)$ and any number $\nu \geq 0$, the graph $G^{\ell^p}_{\varepsilon}(f)$ is a subgraph of $G^{\ell^p}_{(\varepsilon,\nu)}(f)$. 
\end{lemma}

We have the following reduction. 

\begin{theorem}\label{th:graph_collapse}
Let $f \colon X \to X$ be a mapping on a metric space $X$. 
For any $\varepsilon_1 < \varepsilon_2 \in (-\infty,-0] \sqcup [0,\infty)$, the Morse graph $G^{\ell^p}_{(\varepsilon_2,\nu)}(f) = (V_2, D_2(\nu))$ is partially obtained from $G^{\ell^p}_{(\varepsilon_1,\nu)}(f) = (V_1, D_1(\nu))$ by vertex collapses associated to $h_{\varepsilon_2}^{\varepsilon_1}$ {\rm(i.e.} $G^{\ell^p}_{(\varepsilon_1,\nu)}(f) \searrow^{\hspace{-3pt}h_{\varepsilon_2}^{\varepsilon_1}} G^{\ell^p}_{(\varepsilon_2,\nu)}(f)${\rm)}, where the mapping $h_{\varepsilon_2}^{\varepsilon_1} \colon V_1 \to V_2$ is defined by assigning to each $X_1 \in V_1$ the $\varepsilon_2$-chain recurrent component of $\mathrm{CR}^{\ell^p}_{\varepsilon_2}(f)$ containing $X_1$. 
\end{theorem}

\begin{proof}
By definition, we have  $\mathrm{CR}^{\ell^p}_{\varepsilon_1}(f) \subseteq \mathrm{CR}^{\ell^p}_{\varepsilon_2}(f)$. 
Therefore, for any $\varepsilon_1$-chain recurrent component $X_1$ of $\mathrm{CR}^{\ell^p}_{\varepsilon_1}(f)$, there are an element $\varepsilon_2 > \varepsilon_1$ and an $\varepsilon_2$-chain recurrent component $X_2$ of $\mathrm{CR}^{\ell^p}_{\varepsilon_2}(f)$ with $X_1 \subseteq X_2$, because of $\mathrm{CR}^{\ell^p}_{\varepsilon_1}(f) \subseteq \mathrm{CR}^{\ell^p}_{\varepsilon_2}(f)$. 
Therefore, the mapping $h_{\varepsilon_2}^{\varepsilon_1} \colon V_1 \to V_2$ is well-defined.
By definition of directed edge of Morse graph with $(\varepsilon,\nu)$-errors, the map $h_{\varepsilon_2,D(\nu)}^{\varepsilon_1} \colon D_1(\nu) \to D_2(\nu)$ defined by 
\[
h_{\varepsilon_2,D(\nu)}^{\varepsilon_1}(X_1,X'_1) := (h_{\varepsilon_2}^{\varepsilon_1}(X_1),h_{\varepsilon_2}^{\varepsilon_1}(X'_1))
\]
 is well-defined. 
%
\end{proof}

The previous result implies Theorem~\ref{main:02'}. 
Notice that the family $(h_{\varepsilon}^{\varepsilon'})_{\varepsilon' < \varepsilon}$ as in the previous theorem satisfies the cocycle condition. 

%

\subsection{Morse hyper-graphs with $(\varepsilon,\nu)$-errors}

Fix a metric space $X$ with a set $\mathcal{X} = \{ X_i \}_{i \in \Lambda}$ of disjoint subsets.
We introduce Morse hyper-graphs with errors as follows. 

\begin{definition}\label{def:Morse_hyp_graph_map_003}
For any $p \in [1,\infty]$ and any number $\nu \geq 0$, a directed hyper-graph $(V, H(\nu))$ with the vertex set $V := \{ X_i \mid i  \in \Lambda \}$, and with the {\bf directed hyper-edge set} 
\[
H(\nu) := \{ (\{X_\lambda \mid \lambda \in A \}, \{X_\lambda \mid \lambda \in W \}) \mid H^{A}_W(\nu) \neq \emptyset \}
\]
is called the {\bf Morse (directed) hyper-graph} of $\mathcal{X}$ with $\nu$-$\ell^p$-errors, where
\[
H^{A}_{W}(\nu) := \left\{ x \in X \middle| \, x \in \bigcap_{\lambda \in A} \bigcup_{\alpha \in X_\lambda} [\alpha]^{\ell^p}_{\nu+}, \,[x]^{\ell^p}_{\nu+} \cap X_k \neq \emptyset \text{ for any }k \in W \right\}
\]
for any index subsets $A,W \subseteq \Lambda$. 
Then such a graph is denoted by $\bm{H^{\ell^p}_{\mathcal{X}}(\nu)}$. 
\end{definition}

\begin{definition}\label{def:Morse_hyp_graph_map02}
For any $\varepsilon \in (-\infty,-0] \sqcup [0,\infty)$ and any number $\nu \geq 0$, the graph $H^{\ell^p}_{\mathcal{X}}(\nu)$ with $\nu$-$\ell^p$-errors is the {\bf Morse (directed) hyper-graph} of $f$ with $(\varepsilon,\nu)$-$\ell^p$-errors if $\mathcal{X}$ is the set of $\varepsilon$-chain recurrent components of  $\mathrm{CR}^{\ell^p}_{\varepsilon}(f)$. 
Then we denote by $\bm{H^{\ell^p}_{(\varepsilon,\nu)}(f)}$ the Morse hyper-graph of $f$ with $(\varepsilon,\nu)$-$\ell^p$-errors. 
\end{definition}

\begin{definition}\label{def:MG_error_04}
The hyper-graph $\bm{H_{(\varepsilon,\nu)}(f)} := H^{\ell^\infty}_{(\varepsilon,\nu)}(f)$ is called the Morse hyper-graph of $f$ with $(\varepsilon,\nu)$-errors. 
Similarly, the hyper-graph $\bm{H^{\Sigma}_{(\varepsilon,\nu)}(f)} := H^{\ell^1}_{(\varepsilon,\nu)}(f)$ is called the Morse hyper-graph of $f$ with $(\varepsilon,\nu)$-$\sum$errors. 
\end{definition}

Notice that Morse (directed) hyper-graphs can be reduced to Morse (directed) graphs. 
In addition, while Morse (directed) hyper-graphs contain more information for applications such as data analysis, Morse (directed) graphs are a more commonly used concept. 
Therefore, both are expected to have their own advantages.

\subsection{Singular limits}

We have the following singular limit behavior.

\begin{theorem}\label{lem:counter_ex03}
For any $\varepsilon>0$, there is a homeomorphism $f \colon A \to A$ on a subset $A \subset \R^2$ with $G^\Sigma_{(\varepsilon,0)}(f) \subsetneq \bigcap_{\nu > 0} G^\Sigma_{(\varepsilon,\nu)}(f)$. 
\end{theorem}

\begin{proof}
Put $A_0 := \{-1,1\} \times \R$ and let $f_0 \colon A_0 \to A_0$ be the identical mapping. 
Define a sequence $(a_n)_{n \in Z}$ as follows: 
\[
a_n = 
\begin{cases}
-1 + 1/n &  (n \in \Z_{<0}) \\
0 &  (n = 0) \\
1 - 1/n &  (n \in \Z_{>0}) 
\end{cases}
\]
%
Put $A_1 := \{(a_n, 2 \varepsilon n) \mid n \in \Z \}$. 
Set $A := A_0 \sqcup A_1$ and define a homeomorphism $f \colon A \to A$ by $f\vert_{A_0} = f_0$ and $f\vert_{A_1}(a_n,2 \varepsilon n) := (a_{n+1}, 2 \varepsilon(n+1))$. 
Then $\mathrm{CR}_{\varepsilon'}^\Sigma(f) = A_0$ for any $\varepsilon' \in [0,\varepsilon]$. 
Put $X_0 := \{-1\} \times \R$ and $X_1 := \{1\} \times \R$. 
The set of $0$-chain recurrent component of $\mathrm{CR}_{0}^\Sigma(f)$ is $\{ X_0, X_1 \}$.  
Then $D^\Sigma(0) = \emptyset$ and $D^\Sigma(\nu) = \{ (X_0, X_1)\}$ for any $\nu > 0$. 
Therefore, the subgraph $\bigcap_{\nu > 0} G^\Sigma_{(\varepsilon,\nu)}(f)$ consists of two vertices and one directed edge, but the subgraph $G^\Sigma_{(\varepsilon,0)}(f)$ consists of two vertices and no directed edges. 
This implies the assertion. 
\end{proof}

%

We would like to know whether $G_{(\varepsilon,0)}(f) = \bigcap_{\nu > 0} G_{(\varepsilon,\nu)}(f)$ for any homeomorphism $f$ and any $\varepsilon>0$.
The proof of Theorem~\ref{lem:counter_ex03} implies the following observation.

\begin{corollary}\label{corollary:counter_ex}
There are a homeomorphism $f \colon A \to A$ on a subset $A \subset \R^2$ and a flow $v$ on $A$ as follows: 
\[
H^\Sigma_{(\varepsilon,0)}(f) \subsetneq \bigcap_{\nu > 0} H^\Sigma_{(\varepsilon,\nu)}(f)
\]
%
\end{corollary}

\section{Coarse chain recurrence for (semi)flows}

In the following three sections, we demonstrate that results analogous to those obtained for the mapping also hold for semiflows. 
In this section, we introduce $\varepsilon$-$\ell^p$-chain recurrence for semiflows for any $p \in [1,\infty]$. 

Let $v$ be a continuous action $v \colon \K \times X \to X$ on a metric space $(X,d)$, where $\K$ is either $\R$ or $\R_{\geq 0}$. 
Notice that if $\K = \R$ (resp. $\K = \R_{\geq 0}$) then $v$ is a flow (resp. semiflow). 

\subsection{$\varepsilon$-$\ell^p$-chains}

We define the $(\varepsilon,T)$-$\ell^p$-chain recurrence as follows. 

\begin{definition}
For any $p \in [1,\infty]$, any positive integer $n \in \Z_{>0}$, any non-negative number $\varepsilon \geq 0$, and any $T > 0$, a sequence $( x_i )_{i=0}^n$ of points $x_i \in X$ is a {\bf $\bm{(\varepsilon,T)}$-$\bm{\ell^p}$-chain} if there are sequences $(\varepsilon_i)_{i=0}^{n-1}$, $(t_i)_{i=0}^{n-1}$ of non-negative numbers with $\varepsilon \geq \| ( \varepsilon_i )_{i=0}^{n-1} \|_{p}$ and $\min_i \{ t_i \} \geq T$  such that $d(v^{t_i}(x_i), x_{i+1}) = \varepsilon_i$ for any $i \in \{ 0,1, \ldots , n-1 \}$. 
\end{definition}

Notice that any $(\varepsilon,T)$-$\ell^\infty$-chain is an $(\varepsilon,T)$-chain and any $(\varepsilon,T)$-$\ell^1$-chain is an $(\varepsilon,T)$-$\sum$chain. 

\subsection{$\varepsilon$-$\ell^p$-chain recurrence}

For any $p \in [1,\infty]$, any non-negative number $\varepsilon \geq 0$, and any $T > 0$, define $[x]^{\ell^p}_{(\varepsilon,T)}$ as follows: 
\[
\begin{split}
[x]^{\ell^p}_{(\varepsilon,T)} :=&  \{ y \in X \mid x \sim^{\ell^p}_{(\varepsilon,T)} y \} 
\\
=&  \{ y \in X \mid \text{ there is an } (\varepsilon,T) \text{-}\ell^p \text{-chain from }x \text{ to } y \}
\end{split}
\]  
Moreover, we define the following binary relation $\sim^{\ell^p}_\varepsilon$ for any $\varepsilon \geq 0$:  
\[
x \sim^{\ell^p}_{\varepsilon} y \text{ if } x \sim^{\ell^p}_{(\varepsilon,T)} y \text{ for any } T > 0
\]
Define $[x]^{\ell^p}_{\varepsilon}$ as follows: 
\[
\begin{split}
[x]^{\ell^p}_{\varepsilon} := & \{ y \in X \mid x \sim^{\ell^p}_{\varepsilon} y \} 
\\
= & \{ y \in X \mid x \sim^{\ell^p}_{(\varepsilon,T)} y \text{ for any } T > 0 \}
\\
= & \bigcap_{T > 0} [x]^{\ell^p}_{(\varepsilon,T)}
\\
= & \{ y \in X \mid \text{For any } T > 0, \text{ there is an } (\varepsilon,T) \text{-}\ell^p \text{-chain from }x \text{ to } y \}
\end{split}
\] 
We have the following observation. 

\begin{lemma}
For any non-negative numbers $\varepsilon_1 < \varepsilon_2$, we have $[x]_{\varepsilon_1} \subseteq [x]_{\varepsilon_2}$. 
\end{lemma}

Moreover, as the mapping case, we define a binary relation as follows.

\begin{definition}
For any $\varepsilon \geq 0$, we define a binary relation $\sim_{\varepsilon +}$ as follows: 
\[
x \sim_{\varepsilon +} y \text{ if } x \sim_{\varepsilon'} y \text{ for any } \varepsilon' > \varepsilon 
\]
\end{definition}

Define $[x]_{\varepsilon +}$ as follows: 
\[
\begin{split}
& [x]_{\varepsilon +} 
\\
:= &\{ y \in X \mid x \sim_{\varepsilon +} y \}  
\\
= &\{ y \in X \mid x \sim_{\varepsilon'} y \text{ for any } \varepsilon' > \varepsilon \}  
\\
= & \{ y \in X \mid x \sim_{(\varepsilon',T)} y \text{ for any }  \varepsilon' > \varepsilon \text{ and any } T > 0 \}
\\
= & \{ y \in X \mid \text{For any } \varepsilon' > \varepsilon \text{ and any } T > 0, 
\\
& \hspace{80pt} \text{ there is an } (\varepsilon',T) \text{-chain from }x \text{ to } y \}
\end{split}
\]
We have the following observations. 

\begin{lemma}\label{lem:02b}
We have the following statements for any $p \in [1,\infty]$: 
\\
{\rm(1)} For any non-negative number $\varepsilon$, we have the following equality: 
\[
[x]^{\ell^p}_{\varepsilon +} = \bigcap_{\varepsilon' > \varepsilon} [x]^{\ell^p}_{\varepsilon'}  = \bigcap_{\varepsilon' > \varepsilon} \bigcap_{T > 0} [x]^{\ell^p}_{(\varepsilon',T)}
\]
{\rm(2)} For any non-negative numbers $\varepsilon_1 < \varepsilon_2$, we have $[x]^{\ell^p}_{\varepsilon_1 +} \subseteq [x]^{\ell^p}_{\varepsilon_2 +}$. 
\end{lemma}

\begin{lemma}\label{lem:equv_rec_flow}
A point $x \in X$ is chain recurrent if and only if $x \sim^{\ell^\infty}_{0 +} x$. 
Similarly, a point $x \in X$ is chain $\sum$recurrent if and only if $x \sim^{\ell^1}_{0 +} x$. 
\end{lemma}

By replacing the map $f$ with the (semi-)flow $v$ in the definitions ($\varepsilon$-$\ell^p$-chain recurrent point and $\mathrm{CR}^{\ell^p}_{\varepsilon}(f)$ in Definition~\ref{def:cr}; $\varepsilon$-chain recurrent point and $\varepsilon$-$\sum$chain recurrent point and $\mathrm{CR}_{\varepsilon}(f)$ and $\mathrm{CR}^\Sigma_{\varepsilon}(f)$ in Definition~\ref{def:cr2}; and $\varepsilon$-chain recurrent component in Definition~\ref{def:crc_map}), one obtains the corresponding definitions (e.g. $\mathrm{CR}^{\ell^p}_{\varepsilon}(v)$, $\mathrm{CR}_{\varepsilon}(v)$, and $\mathrm{CR}^\Sigma_{\varepsilon}(v)$) for (semi-)flows. 
See Definitions~\ref{def:cr_flow_001}, \ref{def:cr2_flows}, and \ref{def:cr_flow_003} in Appendix for details. 
Moreover, we have the following observation. 

\begin{lemma}
For any non-negative numbers $\varepsilon_1 < \varepsilon_2$, we have $\mathrm{CR}^{\ell^p}_{\varepsilon_1}(v) \subseteq \mathrm{CR}^{\ell^p}_{\varepsilon_2}(v)$. 
\end{lemma}
In contrast to the case of mappings, we observe that $(\mathrm{CR}^{\ell^p}_{\varepsilon}(v))_{\varepsilon \geq 0}$ need not form a filtration.

\begin{lemma}\label{lem:no_filtrations_step}
The flow $v \colon \K \times \R \to \R$ by $v(t,x) = x+t$ has no $\varepsilon$-$\ell^p$-chain recurrent points for any $p \in [1,\infty]$ and any $\varepsilon \in \R$. 
\end{lemma}

We call the inferior $\inf \{ \varepsilon \geq 0 \mid \mathrm{CR}^{\ell^p}_{\varepsilon}(v) = X \}$ the {\bf $\bm{\ell^p}$-circulation cost} of $v$. 
Moreover, the $\ell^\infty$-circulation cost of $v$ is called the {\bf circulation cost} of $v$, and the $\ell^1$-circulation cost of $v$ is called the {\bf circulation $\sum$cost} of $v$. 

%
%

\section{Non-gradient property for (semi)flows}

We have the following variant of the non-gradient property. 
Let $v$ be a continuous action $v \colon \K \times X \to X$ on a metric space $X$, where $\K$ is either $\R$ or $\R_{\geq 0}$. 
%

As above, by replacing the map $f$ with the (semi-)flow $v$ in the definitions ($\sim^{\ell^p}_{-\varepsilon +}$ in Definition~\ref{def:neg_grad_map_001}; $\varepsilon$-$\ell^p$-non-gradient and $\mathrm{CR}^{\ell^p}_{-\varepsilon}(f)$ in Definition~\ref{def:cr02}; $\mathrm{CR}_{-\varepsilon}(f)$ and $-\varepsilon$-$\sum$chain recurrent points in Definition~\ref{def:cr2-}; $-\varepsilon$-chain recurrent component in Definition~\ref{def:cr04}), one obtains the corresponding definitions (e.g. $\mathrm{CR}^{\ell^p}_{-\varepsilon}(v)$, $\mathrm{CR}_{-\varepsilon}(v)$, and $\mathrm{CR}^\Sigma_{-\varepsilon}(v)$) for (semi-)flows. 
See Definitions~\ref{def:sim_neg_flow01}--\ref{def:sim_neg_flow04} in Appendix for details. 
%
We obtain the following inclusions. 

\begin{lemma}
For any $p \in [1,\infty]$, we obtain $\mathrm{CR}^{\ell^p}_{-0}(v) \subseteq \mathrm{CR}^{\ell^p}_{0}(v)$. 
\end{lemma}

\begin{proof}
Recall that $\mathrm{CR}^{\ell^p}_{0}(v)$ is the set of $0$-$\ell^p$-chain recurrent points. 
In other words, we have $\mathrm{CR}^{\ell^p}_{0}(v) = \{ x \in X \mid x \sim^{\ell^p}_{0 +} x\}$. 
Therefore, we have $\mathrm{CR}^{\ell^p}_{-0}(v) = \{ x \in X \mid y \sim^{\ell^p}_{0 +} x \text{ for any }y \in O^+(x) \} \subseteq \{ x \in X \mid x \sim^{\ell^p}_{0 +} x\} = \mathrm{CR}^{\ell^p}_{0}(v)$. 
\end{proof}

%

The continuity and compactness imply the following statements. 

\begin{lemma}
Let $v$ be a {\rm(}semi\,{\rm)}flow on a compact metric space. 
Then $\mathrm{CR}_{-0}(v) = \mathrm{CR}_{0}(v) = \mathrm{CR}_{-0}^\Sigma(v) = \mathrm{CR}_{0}^\Sigma(v) = \mathrm{CR}^{\ell^p}_{-0}(v) = \mathrm{CR}^{\ell^p}_{0}(v)$ for any $p \in [1,\infty]$.
\end{lemma}

The following proof is based on the proof of \cite[Lemma~5]{hurley1995chain}. 

\begin{proof}
Fix any point $x \in \mathrm{CR}(v) = \mathrm{CR}_{0}(v) = \mathrm{CR}^\Sigma_{0}(v) = \mathrm{CR}^\Sigma(v)$. 
Then we have the following observation for any $p \in [1,\infty]$: 
\[
\begin{split}
[x]^{\ell^p}_{0} & = \{ y \in X \mid  \text{For any } T > 0, \text{there is a } (0,T)\text{-}\ell^p \text{-chain from }x \text{ to } y \} 
\\
& = O^+(x)
\end{split}
\]
Fix any point $y \in [x]_{0}^\Sigma = [x]_{0} = O^+(x)$. 
Then there is a positive number $\tau > 0$ such that $v(\tau,x)= y$. 

\begin{claim}
If $y \sim_{0_\Sigma +} x$, then $y \sim^{\ell^p}_{0 +} x$ for any $p \in [1,\infty)$. 
\end{claim}

\begin{proof}
For any $p \in [1,\infty)$ and any non-negative numbers $\varepsilon_0, \varepsilon_1, \ldots , \varepsilon_{n-1}$, since $ (\sum_{i=0}^{n-1} \varepsilon_i)^p \geq \sum_{i=0}^{n-1} \varepsilon_i^p$, we have the following inequality: 
\[
\| ( \varepsilon_i )_{i=0}^{n-1} \|_{1} = \sum_{i=0}^{n-1} \varepsilon_i \geq \left( \sum_{i=0}^{n-1} \varepsilon_i^p \right)^{1/p} = \| ( \varepsilon_i )_{i=0}^{n-1} \|_{p}
\]
This implies the assertion. 
\end{proof}

By the previous claim, it suffices to show that $y \sim_{0_\Sigma +} x$. 
We have the following equivalence relations: 
\[
\begin{split}
x \sim_{0_\Sigma +} y & \Longleftrightarrow x \sim_{\varepsilon_\Sigma '} y \text{ for any } \varepsilon' > 0 
\\
& \Longleftrightarrow x \sim_{(\varepsilon',T)} y \text{ for any } \varepsilon' > 0 \text{ and any } T > 0
\\
& \Longleftrightarrow  \text{For any numbers } \varepsilon', T > 0, 
\\
& \hspace{24pt} \text{there is an } (\varepsilon',T)\text{-}{\textstyle \sum} \text{chain from }x \text{ to } y
\\
\end{split}
\]
It suffices to show that for any $\varepsilon', T > 0$, there is a $(\varepsilon',T)$-$\sum$chain from $x$ to $y$.
Fix any $\varepsilon', T > 0$. 
By $x \in \mathrm{CR}^\Sigma(v)$, there is a $(\varepsilon',T+\tau)$-$\sum$chain $( x_i )_{i=0}^n$ from $x$ to $x$ and there are sequences $(\varepsilon_i)_{i=0}^n$, $(t_i)_{i=0}^n$ of non-negative numbers with $\varepsilon' \geq \sum_{i = 0}^n \varepsilon_i$ and $\min_i \{ t_i \} \geq T+\tau$  such that $d(v^{t_i}(x_i), x_{i+1}) = \varepsilon_i$ for any $i \in \{ 0,1, \ldots , n \}$. 
Define $( y_i )_{i=0}^n$ by $y_0 := y$, $y_i := x_i$ for any $i \in \{1,2, \ldots , n\}$ and $(s_i)_{i=0}^n$ by $s_0 := t_0 - \tau$, $s_i := t_i$ for any $i \in \{1,2, \ldots , n\}$. 
Then the sequence $( y_i )_{i=0}^n$ equipped with sequences $(\varepsilon_i)_{i=0}^n$, $(s_i)_{i=0}^n$ of non-negative numbers satisfies $\varepsilon' \geq \sum_{i = 0}^n \varepsilon_i$ and $\min_i \{ s_i \} \geq T$  such that $d(v^{s_i}(y_i), y_{i+1}) = \varepsilon_i$ for any $i \in \{ 0,1, \ldots , n \}$. 
This means that the sequence $( y_i )_{i=0}^n$ is a $(\varepsilon',T)$-$\sum$chain from $y$ to $x$.
\end{proof}

As the mapping case, we have the following singular limit behavior of $\mathrm{CR}_{-\varepsilon}^\Sigma(v)$.

\begin{theorem}\label{th:counter_ex04+}
For any positive number $\varepsilon>0$, there is a flow $v$ on a circle with $\bigcup_{\varepsilon>0} \mathrm{CR}_{-\varepsilon}^\Sigma(v) \subsetneq \mathrm{CR}_{-0}^\Sigma(v) = \mathrm{CR}_{0}^\Sigma(v)$. 
\end{theorem}

In fact, there is such a flow in \S~\ref{sec:c_ex_001}. 
%
Put $\mathrm{CR}^{\ell^p}_{\infty}(v) := X$ for any $p \in [1,\infty]$. 
Moreover, we can obtain the following filtrations. 

\begin{theorem}\label{th:main03-}
Let $v$ be a continuous action $v \colon \K \times X \to X$ on a metric space $X$, where $\K$ is either $\R$ or $\R_{\geq 0}$. 
For any $p \in [1,\infty]$, the family $(\mathrm{CR}^{\ell^p}_{\varepsilon}(v))_{\varepsilon \in \K \sqcup \{\infty\}}$ is a filtration. 
\end{theorem}


Here, the infinity $\infty$ is the maximal element of the totally ordered set $\K \sqcup \{\infty\}$. 
As mentioned above, the filtrations can be considered as an indicator of the persistence of the circulation of points.
In fact, we have the following observation. 

\begin{lemma}
Let $v$ be a continuous action $v \colon \K \times X \to X$ on a metric space $X$, where $\K$ is either $\R$ or $\R_{\geq 0}$. 
For any $p \in [1,\infty]$, if $\mathrm{CR}^{\ell^p}(v) = X$, then $\mathrm{CR}^{\ell^p}_{-\varepsilon}(v) = X$ for any non-negative number $\varepsilon \geq 0$. 
\end{lemma}


\subsection{$\ell^p$-chain recurrent potential and $\varepsilon$-$\ell^p$-chain recurrent diagram of a {\rm(}semi{\rm)}flow}

As the mapping case, we define the $\varepsilon$-$\ell^p$-chain recurrent diagram of (semi)flows to analyze them as follows. 
Let $v$ be a semiflow (resp. flow) on a metric space $X$ and put $I := [-\infty,-0] \sqcup [0,\infty]$ (resp. $I := [-\infty,\infty]$). 

As above, by replacing the map $f$ with the (semi-)flow $v$ in the definitions ($\ell^p$-chain recurrent potential in Definition~\ref{def:debut}; $\ell^p$-chain recurrent diagram $D^{\ell^p}_{\mathrm{CR}}(f)$ in Definition~\ref{def:crd_map_001}; $\mathrm{CR}_{-\varepsilon}(f)$ and $-\varepsilon$-$\sum$chain recurrent points in Definition~\ref{def:crd_map_002}), one obtains the corresponding definitions (e.g. $\ell^p$-chain recurrent potentials $\tau^{\ell^p}_f(x)$, $\mathrm{CR}^{\ell^p}_{-\varepsilon}(v)$, $\mathrm{CR}_{-\varepsilon}(v)$, and $\mathrm{CR}^\Sigma_{-\varepsilon}(v)$) for (semi-)flows. 
See Definitions~\ref{def:sim_neg_flow01}--\ref{def:sim_neg_flow04} in Appendix for details. 

\subsubsection{Examples}

For instance, consider an attracting flow $v_Z \colon \R \times \R_{\geq 0 } \to \R$ generated by a vector field $Z = -x$. 
Then, for any $p \in [1,\infty]$, the $\ell^p$-chain recurrent potential $\tau^{\ell^p}_f(x)$ and the $\ell^p$-chain recurrent diagram $D^{\ell^p}_{\mathrm{CR}}(v_Z)$ of $v_Z$ satisfy 
\[
\tau^{\ell^p}_f(x) = 
\begin{cases}
\vert x \vert &  (x \neq 0) \\
-\infty &  (x = 0)
\end{cases}
\]
\[
D^{\ell^p}_{\mathrm{CR}}(v_Z) = (\{ 0 \} \times [-\infty, 0)) \sqcup \{ (x, \varepsilon) \mid \varepsilon \geq 0, x \in [-\varepsilon, \varepsilon]\}
\]
and the $\ell^p$-circulation cost of $v_Z$ is $\infty$. 

Moreover, consider a repelling flow $v_Y \colon \R \times \R_{\geq 0 } \to \R$ generated by a vector field $Y = x$. 
Then the $\sum$chain recurrent diagram $D_{\mathrm{CR}}^\Sigma(v_Y)$ and chain recurrent diagram $D_{\mathrm{CR}}(v_Z)$ of $v_Y$ are $\R_{\geq 0} \times \{ 0 \}$, and the circulation $\sum$cost and circulation cost of $v_Y$ are $\infty$. 

\begin{figure}[t]
\begin{center}
\includegraphics[scale=0.315]{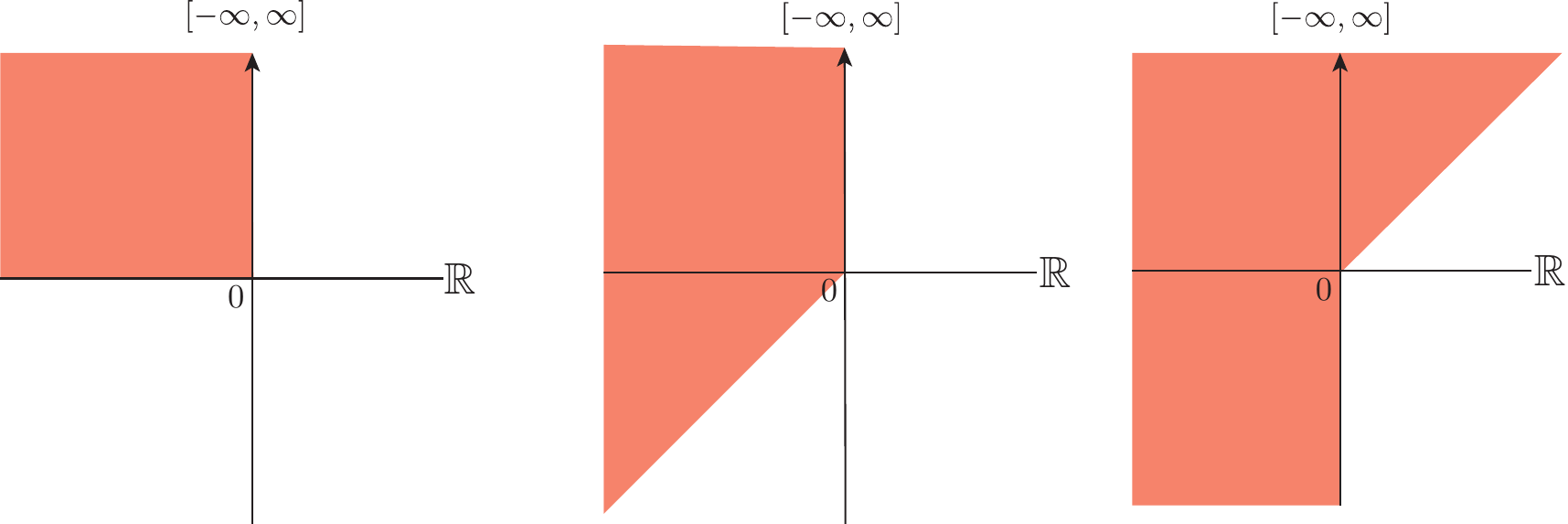}
\end{center} 
\caption{Diagrams $D_{\mathrm{CR}}(v_{\mathrm{rep}})$, $D_{\mathrm{CR}}^\Sigma(v_{\mathrm{rep}})$, and  $D_{\mathrm{CR}}(v_{\mathrm{att}}) = D_{\mathrm{CR}}^\Sigma(v_{\mathrm{att}})$.}
\label{fig:csr_diagram_05}
\end{figure} 

Furthermore, consider flows $v_{\mathrm{rep}}, v_{\mathrm{att}} \colon \R \times \R \to \R$ generated by vector fields $Y_{\mathrm{rep}}, Y_{\mathrm{att}}$, where 
\[
Y_{\mathrm{rep}}(x) = 
\begin{cases}
0 &  (x \leq 0) \\
x &  (x > 0) 
\end{cases}
\]
\[
Y_{\mathrm{att}}(x) = 
\begin{cases}
0 &  (x \leq 0) \\
-x &  (x > 0) 
\end{cases}
\]
are vector fields on $\R$. 
Then the chain ($\sum$)recurrent diagrams of $v_{\mathrm{rep}}$ and $v_{\mathrm{att}}$ are 
\[
D_{\mathrm{CR}}(v_{\mathrm{rep}}) = \R_{\geq 0} \times \R_{\leq 0} \]
\[
\begin{split}
D_{\mathrm{CR}}^\Sigma(v_{\mathrm{rep}}) & = \{ (x, \varepsilon) \mid x \leq \varepsilon < 0 \} \sqcup (\R_{\geq 0} \times \R_{\leq 0}) 
\\
& = \{ (x, \varepsilon) \mid x \leq \varepsilon < 0 \} \sqcup D_{\mathrm{CR}}(v_{\mathrm{rep}})
\end{split}
\]
\[
D_{\mathrm{CR}}(v_{\mathrm{att}}) = D_{\mathrm{CR}}^\Sigma(v_{\mathrm{att}}) = \{ (x, \varepsilon) \mid 0< x \leq \varepsilon \} \sqcup (\R \times \R_{\leq 0})
\]
and the circulation $\sum$costs and circulation costs of $v_{\mathrm{rep}}$ and $v_{\mathrm{att}}$ are $\infty$ as in Figure~\ref{fig:csr_diagram_05}.
The diagrams of these flows suggest that the negative parts of the filtrations $D_{\mathrm{CR}}$ and $D_{\mathrm{CR}}^\Sigma$ correspond to the behaviors near the $\ell^p$-chain recurrent sets.

\section{Morse (hyper-)graphs with finite errors for (semi)flows} 

In this section, we introduce Morse (hyper-)graph with finite errors and Morse (hyper-)graph with finite $\ell^p$-errors for semiflows and flows.  

Let $v$ be a continuous action $v \colon \K \times X \to X$ on a metric space $X$, where $\K$ is either $\R$ or $\R_{\geq 0}$. 
%

\subsection{Morse graphs with $\varepsilon$-$\ell^p$-errors for flows}

As above, by replacing the map $f$ with the (semi-)flow $v$ in the definitions (Morse (directed) graph $G^{\ell^p}_{\varepsilon}(f)$ with $\varepsilon$-$\ell^p$-errors in Definition~\ref{def:MG_error}; $G_{\varepsilon}(f)$ and $G^{\Sigma}_{\varepsilon}(f) $ in Definition~\ref{def:MG_error_02}), one obtains the corresponding definitions (e.g. $G^{\ell^p}_{\varepsilon}(v)$, $G_{\varepsilon}(v)$, and $G^{\Sigma}_{\varepsilon}(v)$) for (semi-)flows. 
See Definitions~\ref{def:MG_error_flow_001} and \ref{def:MG_error_02_flow} in Appendix for details. 
Notice that $G^{\ell^p}_{\varepsilon}(v)$ is a generalization of the Morse graph $G(v)$ of $v$ such that $G^{\ell^\infty}_{0}(v) = G_{0}(v) = G(v)$. 
As in the same argument of the proof of Theorem~\ref{th:graph_collapse_homeomorphism}, we have the following reduction. 

\begin{theorem}\label{th:graph_collapse_flow}
Let $v$ be a flow on a metric space $X$. 
For any $p \in [1,\infty]$ and any $\varepsilon_1 < \varepsilon_2 \in \R$, the Morse graph $G^{\ell^p}_{\varepsilon_2}(v) = (V_2, D_2)$ is partially obtained from $G^{\ell^p}_{\varepsilon_1}(v) = (V_1, D_1)$ by vertex collapses associated to $h_{\varepsilon_2}^{\varepsilon_1}$ {\rm(i.e.} $G^{\ell^p}_{\varepsilon_1}(f) \searrow^{\hspace{-3pt}h_{\varepsilon_2}^{\varepsilon_1}} G^{\ell^p}_{\varepsilon_2}(f)${\rm)}, where the mapping $h_{\varepsilon_2}^{\varepsilon_1} \colon D_1 \to D_2$ is defined by assigning to each $X_1 \in D_1$ the $\varepsilon_2$-chain recurrent component of $\mathrm{CR}^{\ell^p}_{\varepsilon_2}(v)$ containing $X_1$. 
\end{theorem}

%

\subsection{Morse hyper-graphs with $(\varepsilon,\nu)$-errors for (semi)flows}

To define Morse graphs with finite errors for semiflows, we weaken ``connecting'' properties as follows. 
As in the same argument of the proofs of Lemma~\ref{lem:alpha} and Lemma~\ref{lem:omega}, we have the following properties. 

\begin{lemma}\label{lem:alpha_flow}
Let $v$ be a flow on a metric space $X$. 
The following chain of implications holds for any $p \in [1,\infty]$, any \( x \in X \), and any subset \( A \subseteq X \):
\\
{\rm(1)} $\alpha(x) \cap A \neq \emptyset \implies x \in \bigcup_{\alpha \in A} [\alpha]_{0+}^\Sigma $.
\\
{\rm(2)} $x \in \bigcup_{\alpha \in A} [\alpha]_{0+}^\Sigma \implies x \in \bigcup_{\alpha \in A} [\alpha]^{\ell^p}_{0+}$. 
\end{lemma}

\begin{lemma}\label{lem:omega_flow}
Let $v$ be a flow on a metric space $X$. 
The following chain of implications holds for any $p \in [1,\infty]$, any $x \in X$, and any subset $W \subseteq X$:
\\
{\rm(1)} $\omega(x) \cap W \neq \emptyset \implies [x]_{0+}^\Sigma \cap W \neq \emptyset$. 
\\
{\rm(2)} $[x]_{0+}^\Sigma \cap W \neq \emptyset \implies [x]^{\ell^p}_{0+} \cap W \neq \emptyset$.  
\end{lemma}

\subsubsection{Morse graphs with $(\varepsilon,\nu)$-errors}

Fix a set $\mathcal{X} = \{ X_i \}_{i \in \Lambda}$ of disjoint subsets. 

As above, by replacing the map $f$ with the (semi-)flow $v$ in the definitions (Morse graph of $\mathcal{X}$ with $\nu$-errors in Definition~\ref{def:Morse_graph_map_001}; Morse graph $G^{\ell^p}_{(\varepsilon,\nu)}(f)$ with $\ell^p$-$(\varepsilon,\nu)$-errors in Definition~\ref{def:Morse_graph_map02}; Morse graph with $(\varepsilon,\nu)$-($\sum$)errors
 in Definition~\ref{def:MG_error_03}; Morse (directed) hyper-graph $H^{\ell^p}_{\mathcal{X}}(\nu)$ of $\mathcal{X}$ with $\nu$-$\ell^p$-error in Definition~\ref{def:Morse_hyp_graph_map_003};  Morse (directed) hyper-graph $H^{\ell^p}_{(\varepsilon,\nu)}(f)$ with $(\varepsilon,\nu)$-$\ell^p$-errors in Definition~\ref{def:Morse_hyp_graph_map02}; Morse hyper-graph with $(\varepsilon,\nu)$-($\sum$)errors  in Definition~\ref{def:MG_error_04}, one obtains the corresponding definitions (e.g. $G^{\ell^p}_{(\varepsilon,\nu)}(v)$,$G_{(\varepsilon,\nu)}(v)$,  $G^{\Sigma}_{(\varepsilon,\nu)}(v)$, $H^{\ell^p}_{(\varepsilon,\nu)}(v)$, $H_{(\varepsilon,\nu)}(v)$, and $H^{\Sigma}_{(\varepsilon,\nu)}(v)$) for (semi-)flows. 
See Definitions~\ref{def:Morse_graph_flow02-}--\ref{def:MG_error_04_flow} in Appendix for details. 
Lemma~\ref{lem:alpha_flow} and Lemma~\ref{lem:omega_flow} imply the following observation. 

\begin{lemma}
Let $v$ be a flow on a metric space. 
For any $p \in [1,\infty]$, any $\varepsilon \in \R$, and any number $\nu \geq 0$, the graph $G^{\ell^p}_{\varepsilon}(v)$ is a subgraph of $G^{\ell^p}_{(\varepsilon,\nu)}(v)$. 
\end{lemma}

Taking the suspension flows of the homeomorphisms in the proofs of Theorem~\ref{lem:counter_ex03} and Corollary~\ref{corollary:counter_ex} implies the following observation.

\begin{theorem}\label{lem:counter_ex02}
Let $v$ be a flow on a metric space. 
For any $\varepsilon>0$, there is a flow $v$ on a subset on $\R^3$ with $G^\Sigma_{(\varepsilon,0)}(v) \subsetneq \bigcap_{\nu > 0} G^\Sigma_{(\varepsilon,\nu)}(v)$ and $H^\Sigma_{(\varepsilon,0)}(f) \subsetneq \bigcap_{\nu > 0} H^\Sigma_{(\varepsilon,\nu)}(f)$.
%
%
\end{theorem}

\section{Examples}

Roughly speaking, $\varepsilon$-$\sum$ chain recurrence corresponds to control with finite total energy, whereas $\varepsilon$-chain recurrence corresponds to control with finite energy at each step.
First, we observe various phenomena with respect to $\sum$chain recurrence as follows. 

\subsection{Escape from attracting basins of circular flows by controls using finite total energy $\varepsilon$}

The following phenomena serve as toy models of the escape from attracting basins of circular flows by controls that utilize finite total energy.

Define a vector field $X$ on $\R^1$ as follows:
\[
X(x) = 
\begin{cases}
2x \sin \left(\pi/4 - 2/x \right) & x > 0 \\
0 & x = 0 \\
-X(-x) & x < 0 
\end{cases}
\]
Then the set of solutions of the equation $X(x) = 0$ is 
\[
\{ 0 \} \sqcup \left\{ \dfrac{\pm 2}{\pi \left( n+\dfrac{1}{4} \right)} \middle| n \in \Z_{\geq 0} \right\}
\]
because of the odd property of $X$. 
Denote by $v$ the flow generated by $X$ on the line $\R$ such that the subset $\mathrm{CR}_{0}^\Sigma(v) = \mathrm{CR}^\Sigma(v) = \Sv$ is bounded and consists of infinitely many $0$-chain recurrent components.
Moreover, the set difference $\mathrm{CR}_{0}^\Sigma(v) - \{ 0 \} = \{ \pm 2/(\pi(n+1/4)) \mid n \in \Z_{\geq 0} \}$ consists of sinks and sources. 
This implies that $\mathrm{CR}_{\varepsilon}^\Sigma(v)$ for any $\varepsilon > 0$ consists of finitely many $\varepsilon$-chain recurrent components. 
For any $\varepsilon > 8/\pi$, the subset $\mathrm{CR}_{\varepsilon}^\Sigma(v)$ is an interval $[-8/\pi, 8/\pi]$. 

\subsection{Elimination of stagnation of flows by controls using finite total energy $\varepsilon$}

The following phenomenon can be considered a toy model of the circulation of flow by injecting energy to eliminate stagnation.

\subsubsection{Elimination of stagnation of circular flows by controls using finite total energy $\varepsilon$}

Consider a positive number $\mu >1$ and a natural number $k \in \Z_{>0}$. 
Let $p_\pi \colon \R \to \R/2\pi \Z =: \mathbb{S}^1$ be the canonical projection. 
Define an odd function $h \colon \R \to \R$ by $h(x) := x - \mu \sin k x \in T_x \R$. 
Then the derivative $h'(x) = 1 - \mu k \cos k x$ is $2\pi$-periodic and so a vector field $X$ on $\mathbb{S}^1$ defined by $X(p(x)) := dp_{\pi}(x) (- h'(x)) \in T_{p(x)}\mathbb{S}^1$ is well-defined, where $dp_{\pi}(x) \colon T_x \R \to T_{p(x)}\mathbb{S}^1$ is the derivative of $p_\pi$ at $x$. 
Let $v$ be the flow generated by $X$ on the circle $\mathbb{S}^1$ whose singular points consist of $k$ sources and $k$ sinks. 
Put $x_{\mu,k}$ the minimal solution in $(0, \pi/2k]$ of the equation $\cos kx = 1/\mu k$. 
Then $h'(-x_{\mu,k}) = h'(x_{\mu,k}) = 0$ and $x_{\mu,k} + l\pi/k \in \cos^{-1}(1/\mu k)$ for any $l \in Z$. 
Since $h'(0) = 1 - \mu k < 0$, we have $h'(x) < 0$ for any $x \in (-x_{\mu,k}, x_{\mu,k})$. 
By $h(0) = 0$, we obtain $h(x_{\mu,k}) = x_{\mu,k} - \mu \sqrt{1 - (1/\mu k)^2} = x_{\mu,k} -  \sqrt{(\mu k)^2 - 1}/k < 0$ and $h(-x_{\mu,k}) = - h(x_{\mu,k}) > 0$. 
Then $h(-x_{\mu,k}) - h(x_{\mu,k}) = 2 h(-x_{\mu,k}) > 0$.
This implies 
\[
\mathrm{CR}_{\varepsilon}^\Sigma(v) 
\begin{cases}
= p_\pi(\cos^{-1}(1/\mu k)) & \varepsilon = 0 \\
\subsetneq \mathbb{S}^1 & (\varepsilon < 2k h(-x_{\mu,k})) \\
= \mathbb{S}^1  & (\varepsilon \geq 2k h(-x_{\mu,k}))
\end{cases}
\]
because we need energy more than $2 h(-x_{\mu,k})$ to overcome each attracting basin for any point in a small \nbd of a sink that is slightly more positive than the sink.
Moreover, we have that $\mathrm{CR}_{\varepsilon}^\Sigma(v)$ consists of the sink for any small $\varepsilon \ll 0$ in the case $k = 1$, and that $\mathrm{CR}_{\varepsilon}^\Sigma(v) = \emptyset$ for any small $\varepsilon \ll 0$ in the case $k > 1$.

\subsubsection{Elimination of stagnation of toral flows by controls using finite total energy $\varepsilon$}

Let $p \colon \R \to \R/\Z =: \T^1$ be the canonical projection. 
Choose a smooth function $\psi_2 \colon \T^1 \to [0,1]$ with $\psi_2^{-1}(0) = p([0,1/8] \sqcup [1/4,3/8] \sqcup [1/2,5/8])$. 
Consider a vector field $X(x,y) := (1,\psi_2(y))$ on a torus $\T^2 := \R^2/\Z^2$. 
Let $v$ be the flow generated by $X$ on $\T^2$. 
Then $\Pv = p([0,1/8] \sqcup [1/4,3/8] \sqcup [1/2,5/8]) \times \T^1$ and $\mathrm{P}(v) = p((1/8,1/4) \sqcup (3/8,1/2) \sqcup (5/8,1)) \times \T^1 = \T^2 - \Pv$. 
Therefore, we have the following equality: 
\[
\mathrm{CR}_{\varepsilon}^\Sigma(v) = 
\begin{cases}
\emptyset & \left(\varepsilon < -\dfrac{1}{8} \right) \\
p \left( \left[ 0,\dfrac{1}{8} + \varepsilon \right] \sqcup \left[\dfrac{1}{4},\dfrac{3}{8} + \varepsilon \right] \sqcup \left[\dfrac{1}{2},\dfrac{5}{8} + \varepsilon \right]\right) \times \T^1 
& \left( \varepsilon \in \left[-\dfrac{1}{8},0\right] \right) \\
\T^2  & \left(\varepsilon \geq \dfrac{1}{2} \right)
\end{cases}
\]
Note that $\mathrm{CR}_{\varepsilon}^\Sigma(v) \subsetneq \T^2$ for arbitrarily small $\varepsilon > 0$.

%
%

\subsection{Difference between $\varepsilon$-chain recurrence and $\varepsilon$-$\sum$chain recurrence}\label{sec:c_ex_001}

Notice that the $\varepsilon$-$\sum$chain recurrence (resp. circulation $\sum$cost) seems to be more sensitive than the $\varepsilon$-chain recurrence (resp. circulation cost) as follows. 

Let $p \colon \R \to \R/\Z =: \T^1$ be the canonical projection and $\delta >0$ a small positive number. 
Consider a smooth function $\psi \colon \T^1 \to \R_{\geq 0}$ with $\psi^{-1}(0) = p([0,\delta])$. 
Let $v$ be the flow generated by the vector field $\psi \partial/\partial x$ on $\T^1$ whose singular points form $p([0,\delta])$. 
Then $\T^1 = \mathrm{CR}(v) = \mathrm{CR}_{\varepsilon}(v)$ for any $\varepsilon \in \R$. 
Therefore, the chain recurrent diagram of $v$ is $\R \times \T^1$. 
\begin{figure}[t]
\begin{center}
\includegraphics[scale=0.5]{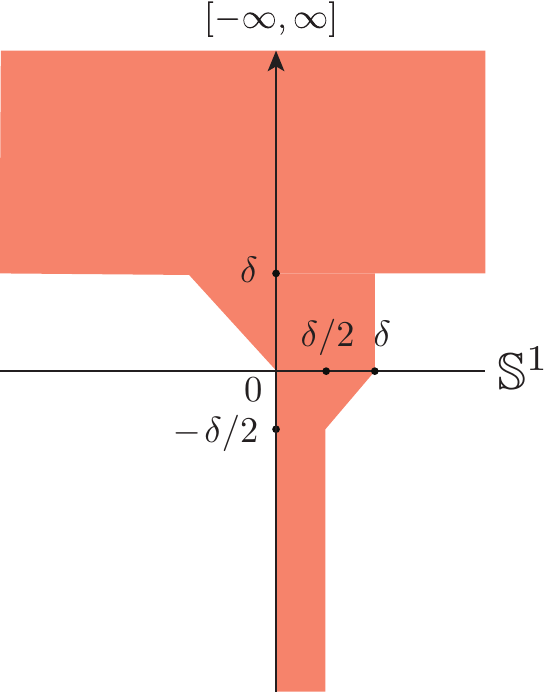}
\end{center} 
\caption{A schematic picture of the $\sum$chain recurrent diagram of a flow $v$.}
\label{fig:csr_diagram_02}
\end{figure} 
On the other hand, we have $\mathrm{CR}_{0}^\Sigma(v) = \mathop{\mathrm{Sing}(v)} = p([0,\delta]) \subsetneq \T^1$. 
Moreover, there is a positive integer $N > 0$ such that 
\[
\begin{split}
0 \in \mathrm{CR}_{-\varepsilon_2}^\Sigma(v) \subsetneq \mathrm{CR}_{-\varepsilon_1}^\Sigma(v)  & \subsetneq  \bigcup_{\varepsilon>0} \mathrm{CR}_{-\varepsilon}^\Sigma(v) = \mathop{\mathrm{Sing}(v)} - \{ \delta \}
\\
& \subsetneq \mathop{\mathrm{Sing}(v)} = \mathrm{CR}_{-0}^\Sigma(v) = \mathrm{CR}_{0}^\Sigma(v)  \\
& \subsetneq \mathrm{CR}_{\varepsilon_1}^\Sigma(v) \subsetneq \mathrm{CR}_{\varepsilon_2}^\Sigma(v)  \subsetneq \T^1
\end{split}
\]
for any numbers $\varepsilon_1 < \varepsilon_2 \in (0,N)$. 
This implies Theorem~\ref{th:counter_ex04+}. 
Furthermore, we have a schematic picture of the $\sum$chain recurrent diagram of $v$ as in Figure~\ref{fig:csr_diagram_02}. 

Let $f$ be the time-one map of $v$ whose fixed points form $p([0,\delta])$. 
Then $\T^1 = \mathrm{CR}(f) = \mathrm{CR}_{\varepsilon}(f)$ for any $\varepsilon \in \R$. 
On the other hand, we have $\mathrm{CR}_{0}^\Sigma(f) = \mathop{\mathrm{Fix}(f)} = p([0,\delta]) \subsetneq \T^1$. 
Moreover, there is a positive integer $N > 0$ such that 
\[
\begin{split}
0 \in \mathrm{CR}_{-\varepsilon_2}^\Sigma(f) \subsetneq \mathrm{CR}_{-\varepsilon_1}^\Sigma(f)  & \subsetneq  \bigcup_{\varepsilon>0} \mathrm{CR}_{-\varepsilon}^\Sigma(f) = \mathop{\mathrm{Fix}(f)} - \{ \delta \}
\\
& \subsetneq \mathop{\mathrm{Fix}(f)} = \mathrm{CR}_{-0}^\Sigma(f) = \mathrm{CR}_{0}^\Sigma(f)  \\
& \subsetneq \mathrm{CR}_{\varepsilon_1}^\Sigma(f) \subsetneq \mathrm{CR}_{\varepsilon_2}^\Sigma(f)  \subsetneq \T^1
\end{split}
\]
for any numbers $\varepsilon_1 < \varepsilon_2 \in (0,N)$. 
This implies Theorem~\ref{th:counter_ex04}. 
Furthermore, we have a schematic picture of the $\sum$chain recurrent diagram of $f$ in Figure~\ref{fig:csr_diagram_01}. 
On the other hand, the chain recurrent diagram of $f$ is the whole space $\R \times \T^1$. 
Moreover, the circulation costs of the circular flow $v$ and the time-one map $f$ are zero and the circulation $\sum$costs of the circular flow $v$ and the map $f$ are $\delta$. 
Thus, the circulation $\sum$costs for the flow and the map represent the energetic thresholds required to transcend stagnation in this example.

\subsubsection{Existence of $\varepsilon$-chains and non-existence of $\varepsilon$-$\sum$chains}\label{sec:ex_e_orbit}

Notice that for any pair of points contained in the same $\varepsilon$-chain recurrent component of the $\varepsilon$-$\sum$chain recurrent set $\mathrm{CR}^{\sum}_{\varepsilon}(f)$ for a mapping $f$, the existence of $\varepsilon$-$\sum$chains between them is generally not necessary.
For instance, for the identical map $1_\R$ on the line $\R$ and any $\varepsilon >0$, there are no $\varepsilon$-$\sum$chains between $a$ and $b$ for any points $a<b \in \R$ with $b-a > \varepsilon$. 

On the other hand, for any pair of points contained in the same $\varepsilon$-chain recurrent component of the $\varepsilon$-chain recurrent set $\mathrm{CR}_{\varepsilon}(f)$ for a mapping $f$, there is an $\varepsilon$-chain between them.

\subsection{Control of flows using finite energy $\varepsilon$ at each step}

We observe various phenomena with respect to chain recurrence as follows. 

\subsubsection{Escape from bounded regions}

 Let $f$ be the time-one map on the line $\R$ of the gradient flow whose height function $H(x) = \cos (x)$. 
Then any $\varepsilon$-chain for a sufficiently small non-negative number $\varepsilon > 0$ is bounded. 
On the other hand, for any point $x \in \R$ and for any $\varepsilon > D$, there is an $\varepsilon$-chain of it is unbounded, where $D := \max \{ \vert f(y) - y \vert \mid y \in \R \} \in (0,\pi)$.
Moreover, we have the following observation:  
\[
\mathrm{CR}_{\varepsilon}(f_R) =
\begin{cases}
\emptyset & \text{if } \varepsilon < - \pi \\
\{ n \pi \mid n \in \Z \} & \text{if } \varepsilon = 0 \\
\mathbb{R} & \text{if } \varepsilon > D.
\end{cases}
\]


\subsubsection{Escape towards positive infinity}

Define a translation $f_R \colon \R \to \R$ by $f(x) := x +R$ for any $R \in \R$. 
Fix any non-positive number $R \leq 0$. 
Then 
\[
\mathrm{CR}_{\varepsilon}(f_R) =
\begin{cases}
\emptyset & \text{if } \varepsilon < -R, \\
\mathbb{R} & \text{if } \varepsilon \geq -R.
\end{cases}
\]
By the observation in \S~\ref{sec:ex_e_orbit}, if $\varepsilon \geq -R$, then we can consider that any point can escape towards positive infinity by using finite energy $\varepsilon$ at each step. 
On the other hand, since any orbit of $f_R$ for any $R >0$ goes to positive infinity, we can consider that any point can escape towards positive infinity without using any energy. 

\subsection{A family of ``non-increasing'' Morse graphs}\label{sec:non-increasing}

\begin{figure}[t]
\begin{center}
\includegraphics[scale=0.5]{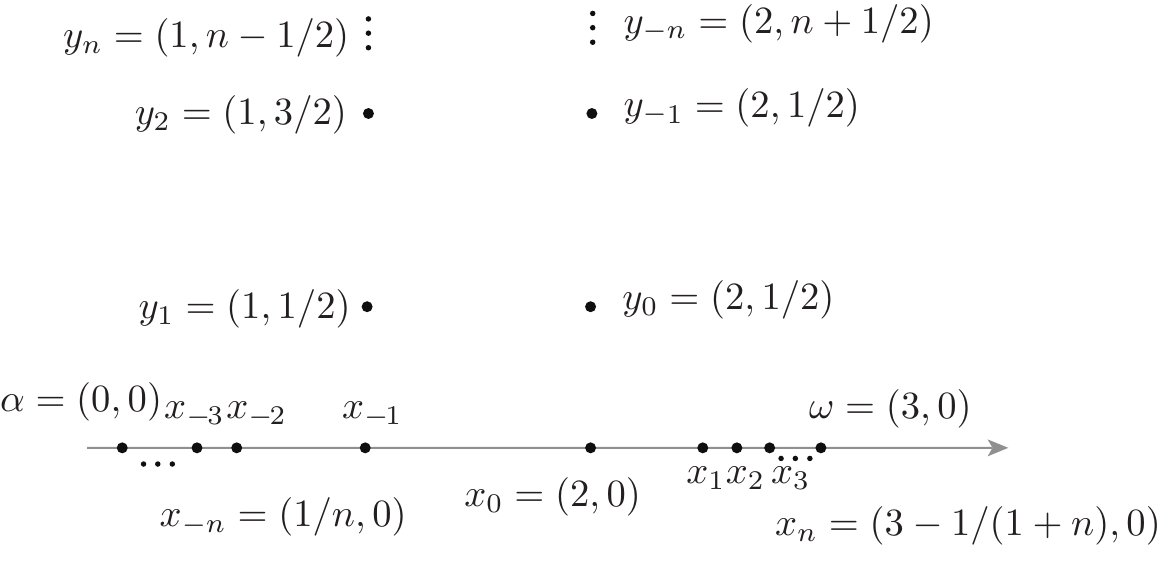}
\end{center} 
\caption{The space $X= \{ \alpha, \omega \} \sqcup \{ x_n, y_n \mid n \in \Z\}$.}
\label{fig:example_001}
\end{figure}

Set $\alpha := (0,0) \in \R^2$ and $\omega := (3,0) \in \R^2$. 
Define $x_{-n} := (1/n,0) \in \R^2$, $x_0 := (2, 0) \in \R^2$, and $x_{n} := (3- 1/(1+n), 0) \in \R^2$ for any $n \in \Z_{>0}$. 
Write $y_{-n} := (2, n+1/2) \in \R^2$, $y_0 := (2,1/2) \in \R^2$, and $y_n := (1,n-1/2) \in \R^2$ for any $n \in \Z_{>0}$. 
Let $X := \{ \alpha, \omega \} \sqcup \{ x_n, y_n \mid n \in \Z\}$ as in Figure~\ref{fig:example_001}.
Define a homeomorphism $f \colon X \to X$ defined by $f(\alpha) = \alpha$, $f(\omega) = \omega$, $f(x_n) = x_{n+1}$, and $f(y_n) = y_{n+1}$. 
Put $X_1 := \{\alpha\}$ and $X_2 := \{\omega\}$. 
Then $\mathrm{CR}_0(f) = \mathrm{CR}(f) = \{ \alpha, \omega \}$ and $G_0'(f) = G'(f) = (V_0,D'_0)$ with $V_0= \{ X_1, X_2 \}$ and $D'_0 = \{(X_1,X_2)\}$. 

On the other hand, put $\varepsilon := 1/2$, and set $X_1(\varepsilon) := X_1 \sqcup \{ x_{-n} \mid -n \leq -2 \}$, $X_2(\varepsilon) := X_2 \sqcup \{ x_{n} \mid n \geq 1 \}$, $X_3(\varepsilon) := \{ x_{-1}, y_0 \}$, $X_4(\varepsilon) := \{x_0 \}$. 
Then $\mathrm{CR}_{\varepsilon}(f) = \bigsqcup_{i=1}^4 X_i(\varepsilon)$ and $G_{\varepsilon}'(f) = (V_{\varepsilon},D_{\varepsilon}')$ with $V_{\varepsilon} = \{ X_i(\varepsilon) \mid i \in \{ 1,2,3,4\} \}$ and $D_{\varepsilon}' = \emptyset$. 
Notice that the mapping $h^0_{\varepsilon} \colon V_0 \to V_{\varepsilon}$ is defined by assigning to each $X \in V_0$ the equivalence class of $\mathrm{CR}_{\varepsilon}(f)$ containing $X$ is well-defined but the map $D'_0 \to D_{\varepsilon}' = \emptyset$ is not well-defined because of $D'_0 \neq \emptyset$.

\section{Final remarks}\label{rem:f_rem}

This paper defines $\varepsilon$-$\ell^p$-chain recurrence for mappings and semiflows. 
On the other hand, one can define variants of these concepts. 
In fact, an infinitesimal version of $\varepsilon$-$\ell^p$-chain recurrence for semiflows can also be defined using the following $\varepsilon$-$\int$chain recurrence. 

\begin{definition}
Let $v$ be a differentiable flow on a manifold $M$. 
For any positive integer $n \in \Z_{>0}$ and any non-negative number $\varepsilon \geq 0$, a differentiable curve $c \colon [0,S] \to M$ is a {\bf $\bm{\varepsilon}$-$\bm{\int}$chain} if the following inequality holds: 
\[
\varepsilon \geq \int_0^S \left\| \left. \dfrac{d v(t,c(s))}{dt} \right\vert_{t = 0} - \dfrac{d c(s)}{ds} \right\| ds 
\]
\end{definition}

As the $\varepsilon$-$\sum$chain recurrence for semiflows, we can define $\varepsilon$-$\int$chain recurrence for semiflows. 
Moreover, a variance of $\varepsilon$-$\ell^p$-chain recurrence for mappings can be defined analogously to $\varepsilon$-$\ell^p$-chain recurrence for flows  as follows:  
Let $f \colon X \to X$ be a mapping. 
To define $(\varepsilon,\infty)$-$\ell^p$-chain recurrence, we define the following concept for mappings. 

\begin{definition}
Let $f$ be a mapping on a metric space $X$ and set a positive integer $N > 0$. 
For any positive integer $n \in \Z_{>0}$ and any non-negative number $\varepsilon \geq 0$, a sequence $( x_i )_{i=0}^n$ of points $x_i \in X$ is a {\bf $(\varepsilon,N)$-$\ell^p$-chain} if there are a sequence $(\varepsilon_i)_{i=0}^{n-1}$ of non-negative integers and a sequence $(t_i)_{i=0}^{n-1}$ of non-negative integers with $\varepsilon \geq \| (\varepsilon_i)_{i=0}^{n-1} \|_p$ and $\min_i \{ t_i \} \geq N$  such that $d(f^{t_i}(x_i), x_{i+1}) = \varepsilon_i$ for any $i \in \{ 0,1, \ldots , n-1 \}$. 
\end{definition}

As the $\varepsilon$-$\ell^p$-chain recurrence for semiflows, we can define $(\varepsilon,\infty)$-$\ell^p$-chain recurrence for mappings as follows. 

\begin{definition}
Let $f$ be a mapping on a metric space $X$. 
A point $x \in X$ is $(\varepsilon,\infty)$-$\ell^p$-chain recurrent if for any integer $N>0$, there is a $(\varepsilon,N)$-$\ell^p$-chain from $x$ to $x$. 
\end{definition}

Note that the $(\varepsilon,\infty)$-$\ell^p$-chain recurrence for mappings has similar properties to $\varepsilon$-$\ell^p$-chain recurrence for semiflows. 
Indeed, for instance, considering a translation $h \colon \R \to \R$ defined by $h(x) := x +1$, though every point of $\R$ is $1$-$\ell^p$-chain recurrent with respect to $h$, no points are $(\varepsilon,\infty)$-$\ell^p$-chain recurrent for any $\varepsilon \in \R$ as Lemma~\ref{lem:no_filtrations_step}. 

Investigating the analysis of which types of phenomena are suitable for what kinds of variants of $\varepsilon$-$\ell^p$-chain recurrence is future work.

\subsection{Remarks on ``gradient-like'' behaviors}

In this paper, we have generalized the concept of recurrence. 
However, the developed framework is also applicable to the analysis of ``gradient-like'' behaviors. 
This broader applicability is exemplified through a practical case study, in which filtrations are constructed from ensemble forecasts of a tropical cyclone \cite{imoto2025filtrations}.
A more thorough investigation of this direction is left for future work.

\subsection{Remarks on practical applications}

In the context of the tropical cyclone study mentioned above, the analysis of such ``gradient-like'' phenomena was conducted using functions analogous to chain recurrent potentials. 
Similarly, analyzing data from various dynamical systems by utilizing chain recurrent diagrams and chain recurrent potentials to enhance the readability of data properties is also left for future work.

\subsection{Generalization to a setting with cost functions and partial mappings}

In this paper, we dealt with mappings on metric spaces. 
However, this setting can be generalized to one that handles partial mappings on sets with a cost function\footnote{A cost function on a set $X$ is a function $c \colon X \times X \to [0,\infty]$ with $c(x,x) = 0$ for any $x \in X$. 
}, allowing for similar discussions. 
In fact, \cite{imoto2025filtrations} deals with partial mappings on sets with a cost function and develops the theory.

\subsection{Two-parameter families}
 
 It is a subject of future work to study two-parameter families of the form
\[
D_{\mathrm{Per}}(f_{\cdot}) := \{ (\mu, \varepsilon, \mathrm{Per}_{\varepsilon}(f_\mu)) \mid \varepsilon \in I \}
\]
for one-parameter families of maps $(f_\mu)_{\mu}$ (e.g. tent maps or unimodal maps). 
In particular, the subspace $\{ (\mu, 0, \mathrm{Per}_{\varepsilon}(f_\mu)) \mid \varepsilon \in I \}$ with $\varepsilon = 0$ corresponds to the classical bifurcation diagram of periodic points, and how these structures change under perturbations with finite error remains to be explored.

\subsection{Application to geodesic flows}

By considering geodesic flows on (especially open) Riemannian manifolds, one can naturally construct filtrations of the manifold associated with recurrence behavior. 
Since geodesic flows can be regarded as models that encompass heat diffusion and other physical transport phenomena, this perspective is expected to connect with a wide range of real-world dynamical behaviors, including symmetry-breaking phenomena. 
Investigating such geometric and applied contexts constitutes another direction for future research.

\subsection{Filtrations and $\ell^p$-chain recurrent diagram on difference equations and differential equations}

We can also introduce $\varepsilon$-$\ell^p$-chain recurrence and related concepts as follows: 
Consider a differential (resp. difference) equation $F(s,x) = 0$ on a manifold $M$ with a cost function $c$ on $M$. 
Then $F^{t}(s, x) \in M$ denotes the position at time $t+s \in \K$ of the solution that was at $x \in M$ at time $s \in \K$, where $\K = \R$ (resp. $\K = \Z$).

\begin{definition}
For any positive integer $n \in \Z_{>0}$ and any non-negative numbers $\varepsilon \geq 0$ and $T \geq 0$, a sequence $( (s_i, x_i) )_{i=0}^n$ of points $(s_i,x_i) \in \K \times M$ is a  {\bf $\bm{(\varepsilon,T)}$-$\bm{\ell^p}$-solution} from $(s_0,x_0)$ to $(s_n,x_n)$ if there are sequences $(\varepsilon_i)_{i=0}^{n-1}$, $(t_i)_{i=0}^{n-1}$ of non-negative numbers with $\varepsilon \geq \| (\varepsilon_i)_{i=0}^{n-1} \|_{p}$, $s_{i+1} = s_i + t_i$, and $\min_i \{ t_i \} \geq T$  such that $c(F^{t_i}(s_i,x_i), x_{i+1}) = \varepsilon_i$ for any $i \in \{ 0,1, \ldots , n-1 \}$. 
\end{definition}

For any $\varepsilon \geq 0$ and $T>0$, define binary relations $\sim^{\ell^p}_{(\varepsilon,T)}$ and $\sim^{\ell^p}_{(\varepsilon+,T)}$
 on $\K \times M$ as follows: 
\[
(s,x) \sim^{\ell^p}_{(\varepsilon,T)} (s',y) \text{ if there is an } (\varepsilon,T) \text{-}\ell^p \text{-solution from }(s,x) \text{ to } (s',y)  
\]
\[
(s,x) \sim^{\ell^p}_{(\varepsilon+,T)} (s',y) \text{ if } (s,x) \sim^{\ell^p}_{(\varepsilon,T)} (s',y) \text{ for any }\varepsilon' > \varepsilon 
\]
We have the following chain recurrence. 

\begin{definition}
For any $\varepsilon \geq 0$, a point $(s,x) \in \K \times M$ is {\bf $\bm{\varepsilon}$-$\bm{\ell^p}$-chain recurrent} if, for any $\varepsilon' > \varepsilon$ and $T>0$, there is a number $s'>s$ with $(s,x) \sim^{\ell^p}_{(\varepsilon',T)} (s',x)$. 
Denote by $\bm{\mathrm{CR}^{\ell^p}_\varepsilon(F)}$ the set of $\varepsilon$-$\ell^p$-chain recurrent points of the differential (resp. difference) equation. 
\end{definition}


Moreover, we define the following binary relation $\sim^{\ell^p}_{\varepsilon}$  for any $\varepsilon \geq 0$:  
\[
(s,x) \sim^{\ell^p}_{\varepsilon} (s',y) \text{ if } (s,x) \sim^{\ell^p}_{(\varepsilon,T)} (s',y) \text{ for any } T > 0
\]
Define $[(s,x)]^{\ell^p}_{\varepsilon}$
 as follows: 
\[
\begin{split}
[(s,x)]^{\ell^p}_{\varepsilon} := &\{ (s',y) \in \K \times M \mid (s,x) \sim^{\ell^p}_{\varepsilon} (s',y) \}
\end{split}
\]

We introduce the following concept, which is ``$\ell^p$-chain recurrence with negative errors'' for differential equations and difference equations. 

\begin{definition}
For any $\varepsilon \geq 0$, we define a binary relation $\sim^{\ell^p}_{-\varepsilon +}$ as follows: 
\[
\begin{split}
(s,x) \sim^{\ell^p}_{-\varepsilon +} (s',y) \text{ if } (s'',z) \sim^{\ell^p}_{\varepsilon' +} (s',y) & \text{ for any } \varepsilon' \in [0,\varepsilon] 
\text{ and any } (s'',z) \in [(s,x)]^{\ell^p}_{\varepsilon'}
\end{split}
\]
\end{definition}


We introduce the following concepts as mappings, each of which is ``$\ell^p$-chain recurrence with negative errors''.

\begin{definition}
For any number $\varepsilon \geq 0$, a point $(s,x) \in \K \times M$ is {\bf $\varepsilon$-non-gradient} if there is a number $s'>s$ with $(s,x) \sim^{\ell^p}_{-\varepsilon +} (s',x)$.
Denote by $\bm{\mathrm{CR}^{\ell^p}_{-\varepsilon}(F)}$ the set of $\varepsilon$-non-gradient points. 
\end{definition}


Put $\bm{\mathrm{CR}^{\ell^p}_{\infty}(F)} := \K \times M$ for any $p \in [1,\infty]$. 
Write $I := (-\infty,-0] \sqcup [0,\infty]$. 
Here, the infinity $\infty$ is the maximal element of the totally ordered set $(-\infty,-0] \sqcup [0,\infty]$. 
Moreover, we can obtain the following filtrations. 

\begin{theorem}\label{th:main03}
For any $p \in [1,\infty]$, the family $(\mathrm{CR}^{\ell^p}_{\varepsilon}(F))_{\varepsilon \in I}$ is a filtration on $\K \times M$. 
\end{theorem}


We introduce the $\ell^p$-chain recurrent diagram of $F$ as follows.

\begin{definition}\label{def:debut_flow_deq}
For any $p \in [1,\infty]$, define the {\bf $\ell^p$-chain recurrent potential} \(\tau^{\ell^p}_F \colon \K \times M \rightarrow [-\infty,-0] \sqcup [0,\infty]\) as follows:
\[
\tau^{\ell^p}_F(s,x) := \inf \{\varepsilon \in I \mid (s,x) \in \mathrm{CR}^{\ell^p}_{\varepsilon}(F) \}
\]
\end{definition}

\begin{definition}
For any $p \in [1,\infty]$, the {\bf $\bm{\ell^p}$-chain recurrent diagram} $D^{\ell^p}_{\mathrm{CR}}(F)$ of $F$ is defined as follows: 
\[
D^{\ell^p}_{\mathrm{CR}}(F) :=\bigsqcup_{\varepsilon \in I} \mathrm{CR}^{\ell^p}_{\varepsilon}(F) \times \{ \varepsilon \} =  \{((s,x),\varepsilon) \in (\K \times M) \times I \mid \tau^{\ell^p}_F(x) \leq \varepsilon \}
\]
\end{definition}


The investigation of the properties of these concepts for differential equations and difference equations is left for future work.

\appendix

\section{Concepts for (semi-)flows}

We state the precise definitions for (semi-)flows as follows. 

\subsection{Coarse chain recurrence for (semi-)flows}

By Lemma~\ref{lem:equv_rec_flow}, we introduce the following variant of recurrence and related concepts. 

\begin{definition}\label{def:cr_flow_001}
For any $p \in [1,\infty]$ and any $\varepsilon \geq 0$, a point $x \in X$ is {\bf $\varepsilon$-$\ell^p$-chain recurrent} if $x \sim^{\ell^p}_{\varepsilon +} x$ {\rm(i.e.} $x \in [x]^{\ell^p}_{\varepsilon +}${\rm)}. 
Denote by $\bm{\mathrm{CR}^{\ell^p}_{\varepsilon}(v)}$ the set of $\varepsilon$-$\ell^p$-chain recurrent points. 
\end{definition}

\begin{definition}\label{def:cr2_flows}
For any $\varepsilon \geq 0$, the subset $\mathrm{CR}^{\ell^\infty}_{\varepsilon}(v)$ is denoted by $\bm{\mathrm{CR}_{\varepsilon}(v)}$ and its points are called {\bf $\bm{\varepsilon}$-chain recurrent points}.
Similarly, the subset $\mathrm{CR}^{\ell^1}_{\varepsilon}(v)$ is denoted by $\bm{\mathrm{CR}^\Sigma_{\varepsilon}(v)}$ and its points are called {\bf $\bm{\varepsilon}$-$\bm{\sum}$chain recurrent points}.
\end{definition}

We define the $\varepsilon$-chain recurrent components to define the Morse (hyper-)graphs with $\varepsilon$-$\ell^p$-errors (see Definitions~\ref{def:Morse_graph_flow02-} and \ref{def:Morse_graph_flow02} for details).  

\begin{definition}\label{def:cr_flow_003}
For any $p \in [1,\infty]$, define a binary relation $\approx_{\mathrm{CR}^{\ell^p}_{\varepsilon+}}$ on $\mathrm{CR}^{\ell^p}_{\varepsilon}(v)$ by $x \approx_{\mathrm{CR}^{\ell^p}_{\varepsilon+}} y$ if $x \sim^{\ell^p}_{\varepsilon+} y$ and $y \sim^{\ell^p}_{\varepsilon+} x$.
The transitive closure of $\approx_{\mathrm{cr}^{\ell^p}_{\varepsilon +}}$ on $\mathrm{CR}^{\ell^p}_{\varepsilon}(v)$ is an equivalence relation on $\mathrm{CR}^{\ell^p}_{\varepsilon}(v)$, denote by $\approx_{\mathrm{CR}^{\ell^p}_{\varepsilon +}}$. 
We call an equivalence class of the chain recurrent set $\mathrm{CR}^{\ell^p}_{\varepsilon}(v)$ with respect to $\approx_{\mathrm{CR}^{\ell^p}_{\varepsilon+}}$ an {\bf $\bm{\varepsilon}$-chain recurrent component} of $\mathrm{CR}^{\ell^p}_{\varepsilon}(v)$. 
\end{definition}

\subsection{Non-gradient concepts for (semi-)flows}

We introduce the following concept, which is ``$\ell^p$-chain recurrence with negative errors'' for the continuous action. 

\begin{definition}\label{def:sim_neg_flow01}
For any $p \in [1,\infty]$ and any $\varepsilon \geq 0$, we define a binary relation $\sim^{\ell^p}_{-\varepsilon +}$ as follows: 
\[
x \sim^{\ell^p}_{-\varepsilon +} y \text{ if } z \sim^{\ell^p}_{\varepsilon' +} y \text{ for any } \varepsilon' \in [0,\varepsilon] \text{ and for any } z \in [x]^{\ell^p}_{\varepsilon'}
\]
\end{definition}

We introduce the following concepts as mappings, each of which is ``$\ell^p$-chain recurrence with negative errors''. 

\begin{definition}\label{def:sim_neg_flow02}
For any $p \in [1,\infty]$ and any number $\varepsilon \geq 0$, a point $x \in X$ is {\bf $\varepsilon$-$\ell^p$-non-gradient} if $x \sim^{\ell^p}_{-\varepsilon +} x$.
Denote by $\bm{\mathrm{CR}^{\ell^p}_{-\varepsilon}(v)}$ the set of $\varepsilon$-$\ell^p$-non-gradient points. 
\end{definition}

\begin{definition}\label{def:cr2-}
For any $\varepsilon \geq 0$, the subset $\mathrm{CR}^{\ell^\infty}_{-\varepsilon}(v)$ is denoted by $\bm{\mathrm{CR}_{-\varepsilon}(v)}$ and its points are called {\bf $\bm{-\varepsilon}$-chain recurrent points}.
Similarly, the subset $\mathrm{CR}^{\ell^1}_{-\varepsilon}(v)$ is denoted by $\bm{\mathrm{CR}^\Sigma_{-\varepsilon}(v)}$ and its points are called {\bf $\bm{-\varepsilon}$-$\bm{\sum}$chain recurrent points}.
\end{definition}

\begin{definition}\label{def:sim_neg_flow04}
For any non-negative number $\varepsilon \geq 0$, define a binary relation $\approx^{\ell^p}_{\mathrm{CR}_{-\varepsilon+}}$ on $\mathrm{CR}^{\ell^p}_{-\varepsilon}(v)$ by $x \approx_{\mathrm{CR}^{\ell^p}_{-\varepsilon+}} y$ if $x \sim^{\ell^p}_{-\varepsilon+} y$ and $y \sim^{\ell^p}_{-\varepsilon+} x$.
The transitive closure of $\approx_{\mathrm{cr}^{\ell^p}_{-\varepsilon +}}$ on $\mathrm{CR}^{\ell^p}_{-\varepsilon}(v)$ is an equivalence relation on $\mathrm{CR}^{\ell^p}_{-\varepsilon}(v)$, denote by $\approx_{\mathrm{CR}^{\ell^p}_{-\varepsilon +}}$. 
We call an equivalence class of the chain recurrent set $\mathrm{CR}^{\ell^p}_{-\varepsilon}(v)$ with respect to $\approx_{\mathrm{CR}^{\ell^p}_{-\varepsilon+}}$ a {\bf $\bm{-\varepsilon}$-chain recurrent component} of $\mathrm{CR}^{\ell^p}_{-\varepsilon}(v)$. 
\end{definition}

\subsection{$\varepsilon$-$\ell^p$-chain recurrent diagram of a {\rm(}semi{\rm)}flow}

We have the following concepts for (semi-)flows. 

\begin{definition}\label{def:debut_flow}
For any $p \in [1,\infty]$, define the {\bf $\ell^p$-chain recurrent potential} \(\tau^{\ell^p}_v \colon X \rightarrow I\) as follows:
\[
\tau^{\ell^p}_v(x) := \inf \{\varepsilon \in I - \{-\infty, \infty \} \mid x \in \mathrm{CR}^{\ell^p}_{\varepsilon}(v) \}
\]
\end{definition}

\begin{definition}\label{def:crd_001}
For any $p \in [1,\infty]$, define the {\bf $\ell^p$-chain recurrent diagram} $D^{\ell^p}_{\mathrm{CR}}(v)$ of $v$ as follows: 
\\
\[
D^{\ell^p}_{\mathrm{CR}}(v) := \bigsqcup_{\varepsilon \in I} \mathrm{CR}^{\ell^p}_{\varepsilon}(v) \times \{ \varepsilon \} = \{(x,\varepsilon) \in X \times I \mid \tau^{\ell^p}_v(x) \leq \varepsilon \}
\]
\end{definition}

\begin{definition}\label{def:crd_002}
The function $\tau^{\ell^\infty}_v$ is denoted by $\bm{\tau_v}$ and is called the {\bf chain recurrent potential}.
The subset $D^{\ell^\infty}_{\mathrm{CR}}(v)$ is denoted by $\bm{D_{\mathrm{CR}}(v)}$ and is called the  {\bf chain recurrent diagram}. 

Similarly, the function $\tau^{\ell^1}_v$ is denoted by $\bm{\tau^\Sigma_v}$ and is called the {\bf $\bm{\sum}$chain recurrent potential}.
The subset $D^{\ell^1}_{\mathrm{CR}}(v)$ is denoted by $\bm{D^\Sigma_{\mathrm{CR}}(v)}$ and is called the  {\bf $\bm{\sum}$chain recurrent diagram}. 
\end{definition}

\subsection{Morse graphs with $\varepsilon$-$\ell^p$-errors for flows}

We introduce the Morse graphs with $\varepsilon$-$\ell^p$-errors for (semi-)flows as follows.

\begin{definition}\label{def:MG_error_flow_001}
For a flow $v$ on a metric space, the Morse graph $G_{\mathcal{X}}(v)$ of $\mathcal{X}$ equipped with the persistent directed edgeset in \S~\ref{sec:Morse_graph} is the {\bf Morse (directed) graph} of $v$ with $\varepsilon$-$\ell^p$-errors if $\mathcal{X}$ is the set of $\varepsilon$-chain recurrent components of $\mathrm{CR}^{\ell^p}_{\varepsilon}(v)$. 
Then we denote by $\bm{G^{\ell^p}_{\varepsilon}(v)}$ the Morse graph of $f$ with $\varepsilon$-$\ell^p$-errors. 
\end{definition}

\begin{definition}\label{def:MG_error_02_flow}
The graph $G_{\varepsilon}(v) := G^{\ell^\infty}_{\varepsilon}(v)$ is called the Morse graph of $v$ with $\varepsilon$-errors. 
Similarly, the graph $G^{\Sigma}_{\varepsilon}(v) := G^{\ell^1}_{\varepsilon}(v)$ is called the Morse graph of $v$ with $\varepsilon$-$\sum$errors. 
\end{definition}

\subsection{Morse graphs with $(\varepsilon,\nu)$-errors}

We introduce Morse graphs with errors as follows. 

\begin{definition}\label{def:Morse_graph_flow02-}
For any $p \in [1,\infty]$ and any number $\nu \geq 0$, a directed graph $(V, D(\nu))$ with the vertex set $V := \{ X_i \mid i  \in \Lambda \}$, and with the {\bf directed edge set} $D(\nu) := \{ (X_j, X_k) \mid D_{j,k}(\nu) \neq \emptyset, j \neq k \in \Lambda \}$ is called the {\bf Morse graph} of $\mathcal{X}$ with $\nu$-errors, where
\[
D_{j,k}(\nu) := \left\{ x \in X  \middle| \, x \in \bigcup_{\alpha \in X_j} [\alpha]^{\ell^p}_{\nu+}, \,[x]^{\ell^p}_{\nu+} \cap X_k \neq \emptyset \right\}.
\]
Then such a graph is denoted by $\bm{G^{\ell^p}_{\mathcal{X},\nu}}$. 
\end{definition}

\begin{definition}\label{def:Morse_graph_flow02}
For any $p \in [1,\infty]$, any $\varepsilon \in \R$, and any number $\nu \geq 0$, the graph $G^{\ell^p}_{\mathcal{X},\nu}$ is the {\bf Morse graph} of $v$ with $(\varepsilon,\nu)$-errors if $\mathcal{X}$ is the set of $\varepsilon$-chain recurrent components of $\mathrm{CR}^{\ell^p}_{\varepsilon}(v)$. 
Then we denote by $\bm{G^{\ell^p}_{(\varepsilon,\nu)}(v)}$ the Morse graph of $v$ with $(\varepsilon,\nu)$-errors. 
\end{definition}

\begin{definition}\label{def:MG_error_03_flow}
For any $p \in [1,\infty]$, the graph $\bm{G_{(\varepsilon,\nu)}(v)} := G^{\ell^\infty}_{(\varepsilon,\nu)}(v)$ is called the Morse graph of $v$ with $(\varepsilon,\nu)$-errors. 

Similarly, the graph $\bm{G^{\Sigma}_{(\varepsilon,\nu)}(v)} := G^{\ell^1}_{(\varepsilon,\nu)}(v)$ is called the Morse graph of $v$ with $(\varepsilon,\nu)$-$\sum$errors. 
\end{definition}

\subsubsection{Morse hyper-graphs with $(\varepsilon,\nu)$-errors}

We introduce Morse hyper-graphs with errors as follows. 

\begin{definition}
For any $p \in [1,\infty]$, any $\varepsilon \in \R$, and any number $\nu \geq 0$, a directed hyper-graph $(V, H(\nu))$ with the vertex set $V := \{ X_i \mid i  \in \Lambda \}$, and with the {\bf directed hyper-edge set} 
\[
H(\nu) := \{ (\{X_\lambda \mid \lambda \in A \}, \{X_\lambda \mid \lambda \in W \}) \mid H^{A}_W(\nu) \neq \emptyset \}
\]
is called the {\bf Morse hyper-graph} of $\mathcal{X}$ with $\nu$-errors, where
\[
H^{A}_{W}(\nu) := \left\{ x \in X \middle| \, x \in \bigcap_{\lambda \in A} \bigcup_{\alpha \in X_\lambda} [\alpha]^{\ell^p}_{\nu+}, \,[x]^{\ell^p}_{\nu+} \cap X_k \neq \emptyset \text{ for any }k \in W \right\}
\]
for any index subsets $A,W \subseteq \Lambda$. 
Then such a graph is denoted by $\bm{H^{\ell^p}_{\mathcal{X}}(\nu)}$. 
\end{definition}

\begin{definition}\label{def:Morse_hyp_graph_flow02}
For any $p \in [1,\infty]$, any $\varepsilon \in \R$, and any number $\nu \geq 0$, the graph $H^{\ell^p}_{\mathcal{X}}(\nu)$ with $\nu$-errors is the {\bf Morse hyper-graph} of $v$ with $(\varepsilon,\nu)$-errors if $\mathcal{X}$ is the set of $\varepsilon$-chain recurrent components of $\mathrm{CR}^{\ell^p}_{\varepsilon}(v)$. 
Then we denote by $\bm{H^{\ell^p}_{(\varepsilon,\nu)}(v)}$ the Morse hyper-graph of $v$ with $(\varepsilon,\nu)$-errors. 
\end{definition}

\begin{definition}\label{def:MG_error_04_flow}
The hyper-graph $\bm{H_{(\varepsilon,\nu)}(v)} := H^{\ell^\infty}_{(\varepsilon,\nu)}(v)$ is called the Morse hyper-graph of $v$ with $(\varepsilon,\nu)$-errors. 

Similarly, the hyper-graph $\bm{H^{\Sigma}_{(\varepsilon,\nu)}(v)} := H^{\ell^1}_{(\varepsilon,\nu)}(v)$ is called the Morse hyper-graph of $v$ with $(\varepsilon,\nu)$-$\sum$errors. 
\end{definition}

\subsection*{Acknowledgments}
The author thanks Yusuke Imoto for their useful comments. 

\bibliographystyle{abbrv}
\bibliography{yt20211011}

\end{document}